\documentclass{amsart}
\usepackage{amssymb}
\usepackage{latexsym}
\usepackage{enumerate}

\usepackage{graphicx}

\usepackage{amsmath}
\usepackage{amsthm}
\usepackage{amssymb}
\usepackage{enumerate}
\usepackage{graphics}
\usepackage{latexsym}

\newcommand\card{\operatorname{card}}

\theoremstyle{plain}
\newtheorem{thm}{Theorem}

\newtheorem{rem}[thm]{Remark}
\newtheorem{prop}[thm]{Proposition}
\newtheorem{lem}[thm]{Lemma}
\newtheorem{res}[thm]{Result}

\numberwithin{equation}{section} \numberwithin{thm}{section}
\pdfoutput=1

\begin{document}

\title[Scattering for $3D$ NLKG below energy norm]{Scattering of rough solutions of the nonlinear Klein-Gordon equations in 3D}

\author{Soonsik Kwon}
\address{Korea Advanced Institute of Science and Technology, 291 Daehak-ro, Yuseong-gu, Daejeon, Korea 34141 }
\email{soonsikk@kaist.edu}

\author{Tristan Roy}
\address{Department of Mathematics, Nagoya University, Furocho, Chikusaku, Nagoya, Japan; Postal Code: 464-8602 }
\email{tristanroy@math.nagoya-u.ac.jp}



\vspace{-0.3in}

\begin{abstract}
We prove scattering of solutions below the energy norm of the nonlinear Klein-Gordon equation in 3D with a defocusing power-type nonlinearity that is
superconformal and energy subcritical: this result extends those obtained in the energy class \cite{brenscat,nakanash1,nakanash2} and those obtained below the energy norm under the additional assumption of spherical symmetry \cite{triroyradnlkg}. In order to do that, we generate an exponential-type decay estimate in $H^{s}$, $s<1$, by means of concentration \cite{bourgjams} and a low-high frequency decomposition \cite{bourgbook,almckstt} : this is the starting point to prove scattering. On low frequencies we modify the arguments in \cite{nakanash1,nakanash2}; on high frequencies we use the smoothing effect of the solutions to control the error terms: this, combined with an almost conservation law, allows to prove this decay estimate.
\end{abstract}

\maketitle

\section{Introduction and Theorem}

In this paper we consider the defocusing nonlinear Klein-Gordon equation on $\mathbb{R}^{3}$:
\begin{equation}
\partial_{tt} u  - \Delta u + u  =  -|u|^{p-1}u
\label{Eqn:NlkgWdat}
\end{equation}
with data $u(0)=u_{0}$, $\partial_{t}u(0)=u_{1}$  lying in $H^{s}$, $H^{s-1}$ respectively. \\
We are interested in the strong solutions of the defocusing nonlinear Klein-Gordon equation on some interval $[0,T]$ i.e maps $u$, $\partial_{t} u$
that lie in $C \left([0, \, T], \, H^{s}(\mathbb{R}^{3}) \right)$, $C \left( [0, \,T], \, H^{s-1} ( \mathbb{R}^{3}) \right)$ respectively and that satisfy
\begin{equation}
u(t)  = \cos{(t \langle D \rangle)} u_{0} + \frac{\sin(t \langle D \rangle)}{\langle D \rangle} u_{1} - \int_{0}^{t} \frac{\sin \left(
(t-t^{'}) \langle D \rangle \right)}{\langle D \rangle} |u|^{p-1}(t^{'}) u(t^{'}) \, dt^{'}.
\label{Eqn:StrongSol}
\end{equation}
The defocusing nonlinear Klein-Gordon equation is closely related to the defocusing nonlinear wave equation:
\begin{equation}
\partial_{tt} v - \triangle v =  - |v|^{p-1} v
\label{Eqn:DefocPWaveEq}
\end{equation}
with data $v(0):=v_{0}$, $\partial_{t} v(0):=v_{1}$. (\ref{Eqn:DefocPWaveEq}) enjoys the following scaling property
\begin{equation}
v(t,x)  \rightarrow \frac{1}{\lambda^{\frac{2}{p-1}}} v \left( \frac{t}{\lambda}, \frac{x}{\lambda} \right), \quad
v_{0}(x)  \rightarrow \frac{1}{\lambda^{\frac{2}{p-1}}} v_{0}  \left( \frac{x}{\lambda} \right), \quad
v_{1}(x)  \rightarrow \frac{1}{\lambda^{\frac{2}{p-1}+1}} v_{1} \left( \frac{x}{\lambda} \right).
\label{Eqn:WaveScaling}
\end{equation}
We define the critical exponent $s_{c}:= \frac{3}{2} - \frac{2}{p-1}$. One can check that the $\dot{H}^{s_{c}} \times \dot{H}^{s_{c}-1}$ norm of
$(u_{0},u_{1})$ is invariant under the scaling transformation (\ref{Eqn:WaveScaling}) \footnote{ Here $\dot{H}^{m}$ denotes the standard homogeneous
Sobolev space endowed with the norm $\| f \|_{\dot{H}^{m}}: = \| D^{m} f \|_{L^{2}(\mathbb{R}^{3})}$}. (\ref{Eqn:NlkgWdat}) is known
to be locally well-posed in $H^{s} \times H^{s-1}$, $s \geq s_c$, $p \geq \frac{7}{3}$ by using an iterative argument.
If $p=5$ then $s_{c}=1$ and we say that the nonlinearity $|u|^{p-1} u$ is $\dot{H}^{1}$ (or energy) critical. If $p=\frac{7}{3}$ then
$s_{c}=0$ and the we say that the nonlinearity is $L^{2}$ (or mass) critical. If $p=3$ then $s_{c}=\frac{1}{2}$ and we say that the
nonlinearity is conformal. If $\frac{7}{3} < p <5$ then we say that the regime is mass supercritical-energy subcritical. If $3<p<5$ then we say that the
regime is superconformal and energy subcritical. \\
It is well-known that smooth solutions of (\ref{Eqn:NlkgWdat}) have a conserved energy
\begin{align}
\nonumber
E(u(t))   := &\frac{1}{2} \int_{\mathbb{R}^{3}} \left| \partial_{t} u (t,x) \right|^{2} \, dx  + \frac{1}{2} \int_{\mathbb{R}^{3}} | \nabla  u
(t,x)|^{2} \, dx + \frac{1}{2} \int_{\mathbb{R}^{3}} | u(t,x) |^{2} \, dx \\
& + \frac{1}{p+1} \int_{\mathbb{R}^{3}} |u(t,x)|^{p+1} \, dx.
\label{Eqn:DefEnergy}
\end{align}
In fact by standard limit arguments the energy conservation law remains true for solutions $(u,\partial_{t}u) \in H^{s} \times H^{s-1}$, $s \geq 1$. Since the lifespan of the local solution depends only on the $H^{s} \times H^{s-1}$ norm of the initial data $(u_{0},u_{1})$ (see \cite{linsog})
for $s>s_{c}$, then it suffices to find an a priori pointwise in time bound in $H^{s} \times H^{s-1}$ of the solution $(u,\partial_{t}u)$ in order to establish global well-posedness. \\
The long-time behavior in the energy space (i.e with data $(u_{0},u_{1}) \in H^{1} \times L^{2}$ ) has attracted much attention from the community.
The energy captures the evolution in time of the $H^{1} \times L^{2}$ norm of the solutions. Since it is conserved we have global existence of
solutions of (\ref{Eqn:NlkgWdat}) in the energy space for all dimension $n$ and for all exponent $p$ that is mass-supercritical and energy-subcritical, i.e
$ 1 + \frac{4}{n} < p < 1 +\frac{4}{n-2} $. The next stage is to understand the asymptotic behavior of the solutions of (\ref{Eqn:NlkgWdat}) in the energy space. The scattering, i.e the linear asymptotic behavior, was proved in \cite{bren,brenscat,ginebvelo,nakanash1, nakanash2,pecher1, pecher2}
for all dimension $n$. \\
The long-time behavior below energy norm (i.e with data in $H^{s} \times H^{s-1}$, $s <1$) has also received much attention from the community. The
global existence of solutions of (\ref{Eqn:NlkgWdat}) has been investigated in \cite{bodaomiao}. The scattering of solutions of (\ref{Eqn:NlkgWdat}) with radial data and in dimension $3$ has been studied in \cite{triroyradnlkg}. More precisely it was proved that the asymptotic behavior for spherical solutions is linear for $3<p<5$ and in $H^{s} \times H^{s-1}$, $ \bar{\bar{s}}:= \bar{\bar{s}}(p) < s < 1$ where
$$  \bar{\bar{s}} :=
\begin{cases}
1- \frac{(5-p)(p-3)}{2(p-1)(p-2)}, \, 4 \geq  p > 3 \\
1- \frac{(5-p)^{2}}{2(p-1)(6-p)}, \, 5 > p \geq 4  . \end{cases}    $$
In this paper we are interested in proving scattering results for general data below the energy norm and in dimension $3$. The main result of
this paper is the following one:

\begin{thm}

Let $5 > p > 3$, $ A \gtrsim 1 $ \footnote{The scattering for small data (i.e $A \ll 1$) is well-known. The proof is also contained in the proof of
our theorem.} and
$(u_{0},u_{1}) \in H^{s} \times H^{s-1}$ such that
\begin{equation}
\| (u_{0},u_{1}) \|_{H^{s} \times H^{s-1}}  \leq A.
\label{Eqn:UpperBdInit}
\end{equation}
Then there exists $\tilde{s}:= \tilde{s}(A,p) < 1$ such that $\tilde{s} \rightarrow 1$ as $A \rightarrow \infty$ and
such that the solution of \eqref{Eqn:NlkgWdat} with data $(u_{0},u_{1})$ exists for all time $T$ and scatters as $T$ goes to infinity, i.e there exists $\left( u_{+,0}, u_{+,1} \right) \in H^{s} \times H^{s-1}$ such that
\begin{equation*}
\lim \limits_{T \rightarrow \infty} \| \left( u(T) , \partial_{t} u(T) \right) - K(T) (u_{+,0}, u_{-,0} ) \|_{H^{s} \times H^{s-1}}  = 0.
\end{equation*}
Here,
\begin{equation}
K(T):= \left(
\begin{array}{ll}
\cos{(T \langle D \rangle)} & \frac{\sin{(T \langle D \rangle)}}{\langle D \rangle} \\
- \langle D \rangle \, \sin{(T \langle D \rangle)} & \cos{ ( T \langle D \rangle)}
\end{array}
\right).
\nonumber
\end{equation}



\label{Thm:Scat}
\end{thm}

\begin{rem}
In fact, as $A \rightarrow \infty$, then there exists a constant $\tilde{\alpha}:= \tilde{\alpha}(p) > 1$ such that one can choose $\tilde{s}$ depending on
$A$ in the following fashion:
\begin{equation}
\tilde{s} = 1 - \frac{1} {\tilde{\alpha}^{\tilde{\alpha}^{ ....^{\tilde{\alpha}}}}}.
\label{Eqn:Formtildes}
\end{equation}
Here the height of the tower is $ \sim A^{\tilde{\alpha}}$.
\end{rem}

\section{Notation}
\label{Sec:Notation}

\subsection{General notation}
\noindent We set some general notation that appear throughout the proof. \\
\\
If $x \in \mathbb{R}$, then $\langle x \rangle:= (1 + |x|^{2})^{\frac{1}{2}}$, $x^{+}$ is a slightly larger number than $x$, $x^{++}$ is a
slightly larger number than $x^{+}$, $x^{-}$ is a slightly smaller number than $x$, and $x^{--}$ is a slightly smaller number than $x^{-}$. \\ \\
Let $f$ be a function differentiable in time and
smooth in space. We write $F(f)$ for the following function

\begin{equation*}
F(f)  := |f|^{p-1} f.
\end{equation*}
Given $J$ a time interval, we denote by $X^{J}_{f}$ the following number

\begin{equation*}
X^{J}_{f}(t):= - \int_{J} \frac{\sin{((t-s)\langle D \rangle)}}{ \langle D \rangle} (f(s)) \, ds
\end{equation*}
\\
\\
We denote by $W$ the set of wave admissible points, i.e
\begin{equation*}
W  := \left\{ (q,r), \, (q,r) \in (2,\infty] \times [2,\infty), \, \frac{1}{q} + \frac{1}{r} \leq \frac{1}{2}  \right\}.
\end{equation*}
Given $m \in [0,1]$, we say that $(q,r)$ is $m-$wave admissible if
$$(q,r) \in W \quad \text{and} \quad
\frac{1}{q} + \frac{3}{r} =\frac{3}{2}-m$$
and $q \geq 2+$ if $m=1$.We denote by $W_{m}$ the set of $m$-wave admissible points. We denote by $\tilde{W}$ the dual set of $W$, i.e
\begin{equation*}
\tilde{W}  := \left\{ (\tilde{x}, \tilde{y}), \, \exists (x,y) \in  W, \,  \frac{1}{x} + \frac{1}{\tilde{x}}=1, \, \frac{1}{y} + \frac{1}{\tilde{y}}=1 \right\}.
\end{equation*}
A graphical representation of these sets is given on Figure \ref{Fig:StrichGraph}.
\subsection{The multiplier $I$, the numbers $Z$, and the mollified energies $E(If)$ }

We introduce the multiplier $I$. \\
\\
The proof of Theorem \ref{Thm:Scat} involves the multiplier $I$ defined as follows:
\begin{equation*}
\widehat{If}(t,\xi):=m(\xi) \hat{f}(t,\xi),
\end{equation*}
where
\begin{itemize}
\item $m(\xi):= \eta \left( \frac{\xi}{N} \right)$
\item $\eta$ is a smooth, radial, nonincreasing function in $|\xi|$ such that
$\eta(\xi):=1$, $|\xi| \leq 1$ and $\eta(\xi):=\frac{1}{|\xi|^{1-s}}$, $|\xi| \geq 2$
\item $N \gg  1$ is a parameter.
\end{itemize}
Throughout the paper we choose $(N,s)$ such that
\begin{equation}
N^{1-s}  \sim 1.
\label{Eqn:ChoiceNs}
\end{equation}
We shall explain in Section \ref{Section:Ideas} why this choice of $(N,s)$ is natural. \\

We introduce some numbers that we constantly use in the proof. Given $J$ a time interval, let
\begin{equation*}
Z_{m,s}(J,f) := \sup_{ \substack {(q,r)-m \\ wave \, adm} } \| \partial_{t} \langle D \rangle^{-m} I  f \|_{L_{t}^{q} L_{x}^{r}(J)} + \|
\langle D \rangle^{1-m} I f \|_{L_{t}^{q} L_{x}^{r}(J)},
\end{equation*}

\begin{equation*}
Z(J,f) := \sup_{m \in [0,1]} Z_{m,s}(J,f),
\end{equation*}
and

\begin{equation}
\begin{array}{ll}
\tilde{Z}(J,f) &  :=
\| \partial_{t} \langle D \rangle^{-\frac{1}{2}} I f \|_{L_{t}^{4} L_{x}^{4} (J)} +
\| \langle D \rangle^{ 1 - \frac{1}{2}} I f \|_{L_{t}^{2(p-1)-} L_{x}^{\frac{6(p-1)}{2p-3} + } (J)} +
\| I f \|_{L_{t}^{2(p-1)-} L_{x}^{2(p-1)+} (J)} \\
& + \| \langle D \rangle^{1- \frac{2}{p}} I f  \|_{L_{t}^{p} L_{x}^{\frac{2p}{p-2}} (J) }
+ \| \langle D \rangle^{1-\frac{3p-5}{2p}} I f  \|_{L_{t}^{p} L_{x}^{2p} (J)} +
\| \langle D \rangle^{1-\frac{1}{2}} I f   \|_{L_{t}^{4+} L_{x}^{4-} (J) } \\
& + \| \langle D \rangle^{1- \frac{p-3}{2} } I f \|_{L_{t}^{\frac{4}{p-3}} L_{x}^{\frac{4}{5-p} }(J)}
+ \| \langle D \rangle^{1-1} I f \|_{L_{t}^{\frac{4(p-1)}{7-p}} L_{x}^{\frac{4(p-1)}{p-3}} (J) }  +
\| \langle D \rangle^{1-s_{c}} I  f \|_{L_{t}^{2(p-1)-} L_{x}^{2(p-1)+} (J)}.
\end{array}
\nonumber
\end{equation}

\begin{rem}
Here $(2(p-1)-,2(p-1)+)$ is a small variation of $(2(p-1),2(p-1))$ such that
$(2(p-1)-,2(p-1)+)$ is $s_{c}$-wave admissible: see the point $A$ on Figure \ref{Fig:StrichGraph}. It is necessary to create this variation to make the proof work. This variation creates other variations of bipoints, such as $(4+,4-)$. For the sake of simplification, we will not describe how these variations are created: the reader is invited to do that himself. It is recommended that the
reader ignores all these variations at first reading.
\end{rem}
\vspace{5mm}
We define the mollified energy of $f$ to be the following:

\begin{align*}
E(If(t)) & := \frac{1}{2} \int_{\mathbb{R}^{3}} |\partial_{t} I f (t,x)|^{2} \, dx + \frac{1}{2} \int_{\mathbb{R}^{3}} |I f(t,x)|^{2} \, dx \\
& + \frac{1}{2} \int_{\mathbb{R}^{3}} |\nabla I f(t,x)|^{2} \, dx + \frac{1}{p+1} \int_{\mathbb{R}^{3}} | I f(t,x)|^{p+1} \, dx.
\end{align*}
We also define

\begin{align*}
E_{c}(If(t)) & := \frac{1}{2} \int_{\mathbb{R}^{3}} |\partial_{t} I f(t,x)|^{2} \, dx + \frac{1}{2} \int_{\mathbb{R}^{3}} |If(t,x)|^{2} \, dx \\
& + \frac{1}{2} \int_{\mathbb{R}^{3}} |\nabla I f (t,x)|^{2} \, dx, \\ \\
E(If(t),B(x_{0},R)) & := \frac{1}{2} \int_{B(x_{0},R)} |\partial_{t} If (t,x)|^{2} \, dx + \frac{1}{2} \int_{B(x_{0},R)} |I f(t,x)|^{2} \, dx    \\
& + \frac{1}{2} \int_{B(x_{0},R)} |\nabla I f(t,x)|^{2} \, dx + \frac{1}{p+1} \int_{B(x_{0},R)} |If(t,x)|^{p+1} \, dx,
\end{align*}
with $B(x_{0},R):= \{ x \in \mathbb{R}^{3}, \, |x-x_{0}| \leq R \}$.\\ \\

\subsection{Leibnitz rules}

\noindent We recall some well-known Leibnitz-type rules. Let $ (r,r_{1},r_{2},\tilde{r}_{1},\tilde{r}_{2},\tilde{r}_{3}) \in (1, \infty)^{6}$ be such that

\begin{equation}
\begin{array}{l}
\frac{1}{r} = \frac{p-1}{r_{1}} + \frac{1}{r_{2}} \, and \\
\frac{1}{r} = \frac{p-2}{\tilde{r}_{1}} + \frac{1}{\tilde{r}_{2}} + \frac{1}{\tilde{r}_{3}}.
\end{array}
\nonumber
\end{equation}
If $\alpha \geq 1-s$ then

\begin{equation}
\| \langle D \rangle^{\alpha} I F(f) \|_{L_{x}^{r}} \lesssim \| f \|^{p-1}_{L_{x}^{r_{1}}}
\| \langle D \rangle^{\alpha} I f \|_{L_{x}^{r_{2}}}.
\nonumber
\end{equation}
If $\alpha \geq 0$ then

\begin{equation}
\begin{array}{ll}
\| \langle D \rangle^{\alpha} (F(f) - F(g)) \|_{L_{x}^{r}} & \lesssim \left( \| f \|^{p-1}_{L_{x}^{r_{1}}} + \| g \|^{p-1}_{L_{x}^{r_{1}}} \right)
\| \langle D \rangle^{\alpha}  (f -g) \|_{L_{x}^{r_{2}}}  \\
&  + \left( \| f \|^{p-2}_{L_{x}^{\tilde{r}_{1}}} + \| g \|^{p-2}_{L_{x}^{\tilde{r}_{1}}} \right) \
 \left( \| \langle D \rangle^{\alpha}  f \|_{L_{x}^{\tilde{r}_{2}}} + \| \langle D \rangle^{\alpha} g \|_{L_{x}^{\tilde{r}_{2}}}    \right)
\| f - g \|_{L_{x}^{\tilde{r_{3}}}}.
\end{array}
\nonumber
\end{equation}
The proof of the first estimate is a simple modification of that of the composition rule (see e.g \cite{taylor}). \\
The second estimate comes from the fundamental theorem of calculus $F(f) = F(g) + \int_{0}^{1} F^{'} (f + t (g-f) ) \cdot (g-f) \, dt$
and the product rule (see e.g \cite{taylor}).

\subsection{The Paley-Littlewood decomposition}

\noindent We constantly use throughout the paper the  Paley-Littlewood decomposition. \\
\\
Let $\phi(\xi)$ be a real, radial, nonincreasing function that is equal to $1$ on the unit ball $\left\{ \xi \in \mathbb{R}^{3}: \, |\xi| \leq 1 \right\}$ and that is supported on $\left\{ \xi \in \mathbb{R}^{3}: \, |\xi| \leq 2 \right\}$. Let $\psi$ denote the function $\psi(\xi)  := \phi(\xi) - \phi(2 \xi)$.
If $(M,M_{1},M_{2}) \in 2^{\mathbb{N}}$ are dyadic numbers such that $M_{2} \geq M_{1}$ we define the Paley-Littlewood operators by \footnote{Here
$\hat{f}$ denotes the Fourier transform in space of $f$}

\begin{align*}
\widehat{P_{M} f}(t,\xi) & := \psi \left( \frac{\xi}{M} \right) \hat{f}(t,\xi), \, M > 1 \\
\widehat{P_{1} f }(t,\xi) & := \phi \left( \xi \right) \hat{f}(t,\xi) \\
\widehat{P_{\leq M} f}(t,\xi) & := \phi \left( \frac{\xi}{M} \right) \hat{f}(t,\xi) \\
\widehat{P_{> M} f}(t,\xi) & := \widehat{f}(t,\xi) - \widehat{P_{\leq M} f}(t,\xi) \\
P_{\ll M } f & : = P_{ \leq \frac{M}{128} }f \\
P_{ \gtrsim  M } f & : = P_{ > \frac{M}{128}}f \\
P_{M_{1} < . \leq M_{2}} f  &  := P_{ \leq M_{2}} f -  P_{\leq  M_{1}} f.
\end{align*}
Then we have
\begin{align*}
f(t) &  = \sum_{M \in 2^{\mathbb{N}}} P_{M} f(t),  \\
f(t) &  = P_{\ll  M}f(t) + P_{\gtrsim M} f(t).
\end{align*}

\subsection{The numbers $\alpha$, $C$ and $c$}

All the relevant constants are denoted by the numbers $\alpha:= \alpha(p) \gtrsim 1 $, $C:=C(p) \gtrsim 1 $,
or $c:= c(p)  \ll 1 $. The constants $\alpha$ are exclusively used in expression involving powers (such as  powers of $A$); otherwise, the
constants $C$ and $c$ are used.\\
\\

\begin{tabular}{|c|c|c|c|c|c|}
  \hline
 Constants $C$ & defined in & Constants $c$ & defined in & Constants $\alpha$ & defined in  \\
  \hline
  $C_{in,1}$ &  (\ref{Eqn:BoundInitMolNrj}) & $c_{1}$   &   (\ref{Eqn:InducStr2}) & $\alpha_{1}$  & (\ref{Eqn:InducStr2})   \\
  $C_{in,2}$ &  (\ref{Eqn:BoundInitMolNrj}) & $c_{2}$    &  (\ref{Eqn:SmallIwBarJPrime}) & $\alpha_{2}$ & (\ref{Eqn:TooLarge}) \\    $C_{1}$    &  (\ref{Eqn:TooLarge}) &  $c_{3}$    &  (\ref{Eqn:SepMolNrj}) &  $\alpha_{3}$ & (\ref{Eqn:SepMolNrj})  \\
  $C_{2}$    &  (\ref{Eqn:FreeDecay}) &   $c_{4}$    &  (\ref{Eqn:Cond4}) & $\alpha_{4}$ & (\ref{Eqn:FreeDecay})  \\
  $C_{3}$  &   (\ref{Eqn:Cond4})      & $c_{5}$  &  (\ref{Eqn:InducStr3})  & $\alpha_{5}$  & (\ref{Eqn:FreeDecay})   \\
  $C_{4}$  &    (\ref{Eqn:Boundk})    & $c_{6}$  &  (\ref{Eqn:ControlSupNrjSmall})    & $\alpha_{6}$  & (\ref{Eqn:Cond4}) \\
  $C_{5}$  &    (\ref{Eqn:SizeK})     & $c_{7}$  &   (\ref{Eqn:SizeK})          & $\alpha_{7}$  & (\ref{Eqn:Boundk}) \\
  $C_{6}$  &    (\ref{Eqn:UpperBdIu}) & $c_{8}$  &   (\ref{Eqn:ConcNrjMol})            & $\alpha_{8}$  & (\ref{Eqn:Cond5Res}) \\
           &                          &  $c_{9}$  &  (\ref{Eqn:BoundZpr})              & $\alpha_{9}$  & (\ref{Eqn:InducStr3}) \\
           &                          &  $c_{10}$  & \text{Lemma} \ref{lem:ShortTimePerturb}   & $\alpha_{10}$  & (\ref{Eqn:ControlSupNrjSmall})   \\
           &                          &          &                                   & $\alpha_{11}$  & (\ref{Eqn:SizeK})  \\
           &                          &          &                                   & $\alpha_{12}$ & (\ref{Eqn:SizeK}) \\
           &                          &          &                                   & $\alpha_{13}$ & (\ref{Eqn:ConcNrjMol}) \\
           &                          &          &                                   & $\alpha_{14}$ & (\ref{Eqn:RepartJj}) \\
           &                          &          &                                   & $\alpha_{15}$ & (\ref{Eqn:UpperBdIu}) \\
           &                          &          &                                   & $\alpha_{16}$ & (\ref{Eqn:BoundZpr}) \\
  \hline
\end{tabular}



\section{Preliminary Results}
\noindent In this section we recall some results that we constantly use throughout the proof of Theorem \ref{Thm:Scat}. \\
\\
Let $J:=[a,b]$ be a time interval.\\
\\
The wave Strichartz estimates (see for example \cite{bren,ginebvelostrich,machahira,nakamura} \footnote{see also \cite{triroyradnlkg}}) can be stated as follows: 
\begin{prop}{\textbf{" Strichartz estimates" }} \label{prop:StrEstLtqLxr}
Assume that $w$ satisfies the following Klein-Gordon equation on $J$
\begin{equation}
\left\{
\begin{array}{ccl}
\partial_{tt} w  - \Delta w + w & = & G \\
w(a,x)& = & w_{0}(x)    \\
\partial_{t} w(a,x) & = & w_{1}(x).
\end{array}
\right.
\label{Eqn:KlGH}
\end{equation}
Then, if $m \in [0,1]$, we have
\begin{align}
\| w \|_{L_{t}^{q} L_{x}^{r} (J)} + &\| \partial_{t} \langle D \rangle^{-1} w \|_{L_{t}^{q} L_{x}^{r} (J) } +  \| w \|_{L_{t}^{\infty}
H^{m} (J) } +  \| \partial_{t} w \|_{L_{t}^{\infty} H^{m-1} (J)} \nonumber \\
& \lesssim   \| w_{0} \|_{H^{m}} + \| w_{1} \|_{H^{m-1}} + \| G \|_{L_{t}^{\tilde{q}} L_{x}^{\tilde{r}} (J)},
\label{Eqn:StrNlkg}
\end{align}
under the assumptions
$$ (q,r) \in W_{m}, \qquad (\tilde{q},\tilde{r}) \in \tilde{W}, \qquad \text{and} \quad  \frac{1}{\tilde{q}}+ \frac{3}{\tilde{r}} -2  = \frac{1}{q} + \frac{3}{r}.$$
\end{prop}
\noindent The next proposition shows that the mollified energy at time zero of the solution $u$ of (\ref{Eqn:NlkgWdat}) with data $(u_{0},u_{1})$ satisfying
(\ref{Eqn:UpperBdInit}) is bounded:

\begin{prop}{``\textbf{ Boundedness of mollified energy at time $0$ }'' \cite{triroyradnlkg}} \label{Prop:EstInitMolNrj}
There exist two constants $C_{in,1}$ and $C_{in,2}$ such that \footnote{ Notice that in (\ref{Eqn:BoundInitMolNrj}) we have deliberately chosen to keep the
term $N^{1-s}$. Indeed, we will use (\ref{Eqn:BoundInitMolNrj}) in Section \ref{Section:Ideas} to explain why it is natural to choose
$N^{1-s} \sim 1$ }

\begin{equation}
E(Iu(0)) \leq C_{in,1} N^{2(1-s)} A^{p+1} \leq C_{in,2} A^{p+1}
\label{Eqn:BoundInitMolNrj}
\end{equation}
\end{prop}
$ $\\
The next proposition shows that the variation of a solution of (\ref{Eqn:NlkgWdat}) on
a time interval can be estimated. More precisely
\begin{prop}{\textbf{"Almost Conservation Law "} \cite{triroyradnlkg}} \label{Prop:Acl}
Let $w$ be a solution of (\ref{Eqn:NlkgWdat}). Let $t_{0} \in J$.  Then
\begin{equation}
\left| \sup_{t \in J} E(Iw) - E(Iw(t_{0})) \right|  \lesssim  \int_{J} \int_{\mathbb{R}^{3}}  | \partial_{t} Iw |
| IF(w)-F(Iw) |  \, dx \, dt \\
 \lesssim \frac{Z^{p+1}(J,w)}{N^{\frac{5-p}{2}-}}.
\label{Eqn:Acl}
\end{equation}
\end{prop}
$ $\\
The last proposition allows to control a weighted norm of $w$ on a time interval; more precisely 

\begin{prop}{\textbf{"Almost Morawetz-Strauss estimate"} \cite{triroyradnlkg}} \label{Prop:AlmMor}
Let $w$ be a solution of (\ref{Eqn:NlkgWdat}). Let $\tilde{x} \in \mathbb{R}^{3}$. Then
\begin{equation}
\int_{J} \int_{\mathbb{R}^{3}} \frac{|Iw|^{p+1}}{|x-\tilde{x}|} \, dx dt  \lesssim \sup_{t \in J} E(Iw(t)) + R_{1}(J,w) +
R_{2}(J,w) .
\label{Eqn:AlmMorStrEst}
\end{equation}
with
\begin{align*}
R_{1}(J,w) &: = \int_{J} \int_{\mathbb{R}^{3}} \frac{\nabla Iw.(x - \tilde{x})}{|x - \tilde{x}|} \left( F(Iw) - I F(w) \right) \, dx dt
\\
\text{and } \quad & R_{2}(J,w) := \int_{J} \int_{\mathbb{R}^{3}} \frac{Iw}{|x - \tilde{x}|} \left( F(Iw) - I F(w) \right) \, dx dt.
\end{align*}
Moreover for $i=1,2$ we have

\begin{equation}
| R_{i}(J,w) |  \lesssim  \frac{Z^{p+1}(J,w)}{N^{\frac{5-p}{2}-}}.
\label{Eqn:BdRi}
\end{equation}

\end{prop}

\section{ Ideas}
\label{Section:Ideas}
\noindent In this section we explain the main ideas of this proof. \\
\\
In \cite{nakanash1,nakanash2}, Nakanishi found for some $(q:=q(p),r:=r(p))$ an upper bound of the $L_{t}^{q} L_{x}^{r}$ norm of the solution of
(\ref{Eqn:NlkgWdat}) with data in the energy class by a tower-exponential type bound of the energy of the form
\begin{equation}
 \| u \|_{L_{t}^{q} L_{x}^{r}(\mathbb{R})}  \lesssim  E^{E^{....^{E}}},
\label{Eqn:TowerTypeBoundE}
\end{equation}
where the height of the tower also depends on the energy. This decay estimate was the preliminary step to prove scattering. A natural question is: is it possible to prove decay estimates of this form (or a modified form) for rougher data? If this is possible, then it might help us to prove scattering of solutions of (\ref{Eqn:NlkgWdat}), by analogy with the scattering theory for data in the energy space. This paper gives a positive answer to this
question, at least for $s$ close enough to one. \\
\\
Of course, one cannot use the energy conservation law because the energy can be infinite on $H^{s} \times H^{s-1}$, $s<1$. Instead we introduce the
multiplier $I$ and we work with the mollified energy of $u$ that is finite in these rough spaces and that is
almost conserved: this is the $I$-method (see e.g \cite{almckstt,keeltaoI}), inspired by the
\textit{Fourier truncation method}, designed in \cite{bourgbook}. We aim at proving a decay estimate that is finite in $H^{s} \times H^{s-1}$. Therefore, by analogy with the energy conservation law, our decay estimate should not only depend on $u$ but also on $I$. It was proved in \cite{triroyradnlkg} that, under the additional assumption of radial symmetry, we can control pretty easily the norm $\| Iu \|_{L_{t}^{p+2} L_{x}^{p+2} (\mathbb{R})}$ by combining the ``Almost
Morawetz-Strauss estimate'' (\ref{Eqn:AlmMorStrEst}) with a pointwise decay estimate, namely a radial Sobolev inequality. Unfortunately, such a
pointwise decay estimate does not exist for general data and we shall establish, by means of concentration
\cite{bourgjams}, a tower-exponential bound of the norm $\| I u \|_{L_{t}^{2(p-1)-} L_{x}^{2(p-1)+} (\mathbb{R}) }$. We shall
denote this norm by the \textit{target norm}. \\
\\
In order to prove this bound, the idea is to use the $H^{1}$ theory for frequencies smaller or equal to the parameter $N$ and to control for frequencies larger than $N$ all the errors that appear in the process of generating this estimate. The success of the $I$ method depends on one condition: by choosing
appropriately a parameter $N >> 1$, one can make the variation of the mollified energy on an arbitrarily long time interval
small compare with its initial size at time zero. Assuming that this condition is satisfied for a moment we can neglect the variation of
the mollified energy and we expect to have, by analogy with the bound (\ref{Eqn:TowerTypeBoundE}) of the $L_{t}^{q} L_{x}^{r}$ norm of the
solutions of (\ref{Eqn:NlkgWdat}) with data in the energy class, a tower-exponential bound of the form
\begin{equation}
\| I u \|_{L_{t}^{2(p-1)-} L_{x}^{2(p-1)+} (\mathbb{R})} \lesssim (C_{in,1} N^{2(1-s)} A^{p+1})^{ ...^{C_{in,1} N^{2(1-s)} A^{p+1}}},
\label{Eqn:TowerTypeBoundEIu}
\end{equation}
where we substitute the energy $E$ in (\ref{Eqn:TowerTypeBoundE}) for the initial size of the mollified energy, i.e $C_{in,1} N^{2(1-s)} A^{p+1}$ by
Proposition \ref{Prop:EstInitMolNrj}. The variation of the mollified energy is then estimated by iteration of the almost conservation law (see Proposition \ref{Prop:Acl}) on intervals such that the target norm is small (see Proposition \ref{prop:LocalBd}), using (\ref{Eqn:TowerTypeBoundEIu}). One gets, roughly speaking, a variation of the form
\begin{equation*}
Var  := \frac{ (C_{in,1} N^{2(1-s)} A^{p+1})^{ ...^{C_{in,1} N^{2(1-s)} A^{p+1}}}}{N^{ \frac{5-p}{2}- }}.
\end{equation*}
In order to satisfy our condition, one must make $Var$ small compare with the initial size of the mollified energy: it is natural to choose $(N,s)$ such that
\begin{equation*}
N^{1-s}  \sim 1,
\end{equation*}
and $N:= N(A) $ a large number depending on $A$. \\
\\
Theorem \ref{Thm:Scat} is proved in Section \ref{Sec:ProofThm}. The proof is based upon a modification of the method of induction on
levels of the conserved energy for data in the energy space that is designed in \cite{bourgjams}. Indeed, since the mollified energy is not conserved, we have to modify significantly this method. In particular, we design a relation that allows to control not only the target norm but also the mollified energy of a solution of (\ref{Eqn:NlkgWdat}), assuming that we control its mollified energy at one time: see the definition of $\mathcal{P}(l)$. Then we prove that this relation is true for large levels of mollified energy at one time by induction, using the small mollified energy at one time theory (see Proposition \ref{prop:EstSmallNrj}). Also, we have to make sure that we can make the mollified energy decrease at one time at a non-decreasing rate, in order to reach the small mollified energy level and apply the small mollified energy theory: this is done by introducing the parameters $c_{6}$ and $\alpha_{10}$ in order to make the variation of the mollified energy small enough. The proof of Theorem  \ref{Thm:Scat} relies upon some propositions that we prove in the other sections. In Section \ref{Sec:LocalBd}, we prove some local bounds: these bounds, combined with Proposition \ref{Prop:Acl} (resp.
Proposition \ref{Prop:AlmMor}), allow to estimate by iteration the variation of the mollified energy (resp. an ``Almost Morawetz-Strauss estimate'') on
an arbitrarily long-time interval. In Section \ref{Sec:SepLocMol} and Section \ref{Sec:ProofPerturb} we modify arguments
of \cite{bourgjams,nakanash1,nakanash2} to separate the localized mollified energy and prove a perturbation result. Notice that throughout these sections,
the multiplier $I$ does not commute with the nonlinearity, and one has to prove some commutator estimates, i.e estimates involving the commutator
$I(F(f)) - F(If)$: these estimates are proved in Appendix $A$.

\begin{figure}[!t]
\centering
\includegraphics[height=15cm]{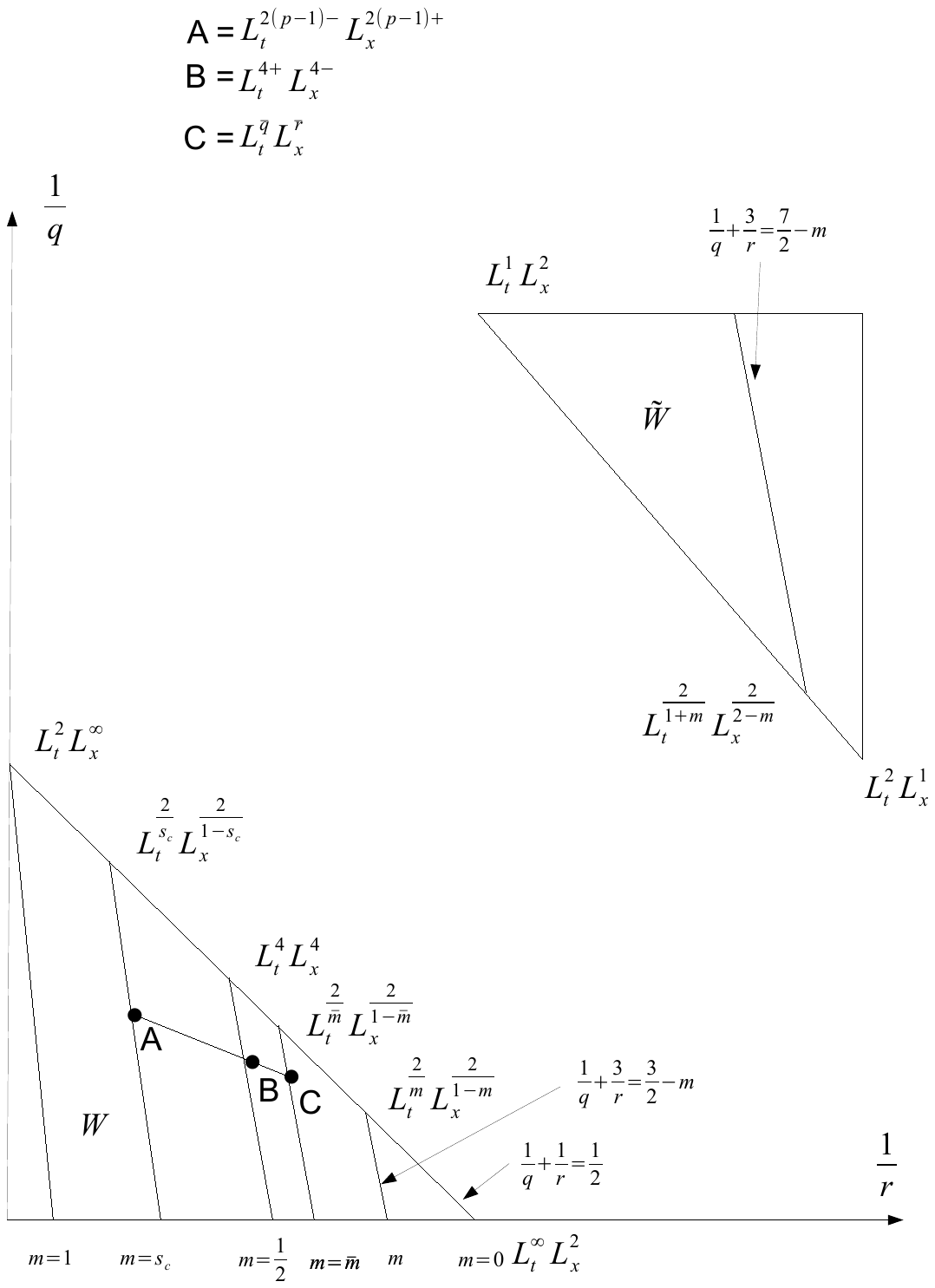}
\caption{The Strichartz graph}
\label{Fig:StrichGraph}
\end{figure}

\section{Proof of Theorem} 
\label{Sec:ProofThm}
\noindent In this section, we prove Theorem \ref{Thm:Scat}.\\
The proof of Theorem \ref{Thm:Scat} relies upon several propositions such as local boundedness, separation of the localized mollified energy, perturbation argument, and small mollified energy theory. We shall prove these propositions in Section \ref{Sec:LocalBd}, \ref{Sec:SepLocMol},
\ref{Sec:ProofPerturb}, and \ref{Sec:EstSmallNrj}.
\subsection{Propositions}
$ $ \\
We consider $J$, $J'$, and $\bar{J}^{'}$ three intervals

\begin{align*}
J:=[a,b] \subset [0, \infty); \, \bar{J}^{'} \subset J' \subset J \cdot
\end{align*}
\\
Let $w$ be a solution of (\ref{Eqn:NlkgWdat}). \\
\\
We assume that there exist $L_{w}$, $X_{w}$ such that

\begin{align}
\sup_{t \in J} E(I w(t))  \leq X_{w} \lesssim A^{p+1},
\label{Eqn:InducNrjw}
\end{align}
and

\begin{align}
\| I w \|_{L_{t}^{2(p-1)-} L_{x}^{2(p-1)+}(J)} \leq L_{w} < \infty.
\label{Eqn:InducStr1}
\end{align}
We will prove in Subsection \ref{Subsec:TheProofThm} that these assumptions are always true.\\
\\
Under these assumptions, one can find constants $\alpha_{1}$ and $c_{1}$ such that if
\begin{align}
\frac{ \langle L_{w} \rangle^{\alpha_{1}} A^{\alpha_{1}}}{ N^{(1-s_{c})-}} \leq c_{1},
\label{Eqn:InducStr2}
\end{align}
then three propositions hold.\\
The first proposition shows that if $\bar{J}^{'}$ is small in the sense of
(\ref{Eqn:SmallIwBarJPrime}), then one can control several norms on this subinterval:
\begin{prop}{\textbf{``Local boundedness''}} \label{prop:LocalBd}
There exists a constant $c_{2}$ such that if
\begin{equation} \label{Eqn:SmallIwBarJPrime}
\begin{array}{ll}
\| I w \|_{L_{t}^{2(p-1)-} L_{x}^{2(p-1)+} (\bar{J}^{'})}   \leq c_{2},
\end{array}
\end{equation}
then
\begin{equation}
\begin{array}{ll}
Z(\bar{J}^{'},w) \lesssim X^{\frac{1}{2}}_{w}.
\end{array}
\label{Eqn:BoundZ}
\end{equation}
\end{prop}
The second proposition shows that if a subinterval is too large in the sense of (\ref{Eqn:TooLarge}), then one can separate the
mollified energy into into two parts: one that is carried by a free Klein-Gordon solution and the other one that
is carried by another solution $w'$ of (\ref{Eqn:NlkgWdat}). The
proof of Proposition \ref{Prop:SepLocMol} relies upon the combination of the $I$-method with a
modification of arguments from Bourgain \cite{bourgjams}, or, more closely, Nakanishi \cite{nakanash1, nakanash2}.
\begin{prop}{\textbf{`` Separation of the localized
mollified energy''}} \label{Prop:SepLocMol}
Let $M \gg 1$. There exist $T \in J^{'}$, and $v$, a solution of the free Klein-Gordon equation, and constants
$c_{3}$, $C_{1}$, $C_{2}$, $\alpha_{2}$,..., $\alpha_{5}$ such that if

\begin{equation}
\| I w \|_{L_{t}^{2(p-1)-} L_{x}^{2(p-1)+} (J^{'})} \geq
C_{1}^{( C_{1}^{A^{\alpha_{2}}} M )^{ C_{1} A^{\alpha_{2}}}},
\label{Eqn:TooLarge}
\end{equation}
then

\begin{align}
E_{c}(Iv) & \lesssim 1,
\label{Eqn:ConservkinetIv} \\
E(Iw'(T)) & \leq  \sup_{t \in J^{'}} E(Iw(t)) -  c_{3} A^{-\alpha_{3}},
\label{Eqn:SepMolNrj} \\
\text{and} \quad \| I v \|_{L_{t}^{2(p-1)-} L_{x}^{2(p-1)+}(J'')} & \leq C_{2} \frac{A^{\alpha_{4}}}{M^{\alpha_{5}}},
\quad \text{for} \quad J'':= (T,b) \quad \text{or} \quad \text{for} \quad  J'':=(a,T).
\label{Eqn:FreeDecay}
\end{align}
Here $w'$ denotes the solution of (\ref{Eqn:NlkgWdat}) such that $w'(T):= w(T) - v(T)$.
\end{prop}
The third proposition shows that if the target norm of $v$ is small in the sense of
(\ref{Eqn:Boundk}), then the target norm of $w$ can be estimated
from that of $w'$ :
\begin{prop}{\textbf{``Perturbation argument''}}
Let $T$, $v$, and $w'$ be defined in the previous proposition. Let $J'':= (T,b)$ (or $J'':=(a,T)$). Assume that $w'$ satisfies (\ref{Eqn:InducNrjw}) and (\ref{Eqn:InducStr1}) (with $w$ substituted for $w'$). Then there exist $c_{4}$, $C_{3}$, $C_{4}$, $\alpha_{6}$, $\alpha_{7}$, $\alpha_{8}$, and $k$ such that if

\begin{equation}
\frac{ C_{3}^{ (A \langle L_{w'} \rangle)^{\alpha_{6}}} (A \langle L_{w} \rangle)^{\alpha_{6}} } {N^{(1-s_{c})-}}  \leq c_{4},
\label{Eqn:Cond4}
\end{equation}

\begin{equation}
k  \leq \frac{1}{C_{4}^{ (A \langle L_{w'} \rangle)^{\alpha_{7}}}} \quad  \text{and} \quad \| I v \|_{L_{t}^{2(p-1)-} L_{x}^{2(p-1)+} (J'')} \leq k,
\label{Eqn:Boundk}
\end{equation}
then

\begin{equation}
\| I  w  \|_{L_{t}^{2(p-1)-} L_{x}^{2(p-1)+} (J'')} \lesssim \left( \langle L_{w'} \rangle A \right)^{\alpha_{8}}.
\label{Eqn:Cond5Res}
\end{equation}
\label{prop:PerturbArg}
\end{prop}
Next we show that if the mollified energy is small enough at one time, then we can have a very good control of
the mollified energy and our target norm:
\begin{prop}{\textbf{`` Small mollified energy theory''}} \label{prop:EstSmallNrj}

Assume that there exists $\tilde{t} \in \mathbb{R}$ such that

\begin{equation}
E(Iw (\tilde{t}))  \ll 1.
\label{Eqn:SmallMolNrjTime}
\end{equation}
Then one can find constants $\alpha_{9}$ and $c_{5}$ such that if

\begin{align}
\frac{ A^{\alpha_{9}}}{ N^{(1-s_{c})-}} \leq c_{5},
\label{Eqn:InducStr3}
\end{align}
then

\begin{equation}
\| I w \|_{L_{t}^{2(p-1)-} L_{x}^{2(p-1)+} (\mathbb{R})}  \lesssim 1,
\label{Eqn:ControlIw}
\end{equation}
and there exist $c_{6}$ and $\alpha_{10}$ such that
\begin{equation}
\sup_{t \in \mathbb{R}} E(Iw(t)) \leq E(Iw(\tilde{t}))(1+
c_{6} A^{- \alpha_{10}}).
\label{Eqn:ControlSupNrjSmall}
\end{equation}
In fact, one can choose $c_{6}$ (resp. $\alpha_{10}$) arbitrarily small (resp. arbitrarily large) in
(\ref{Eqn:ControlSupNrjSmall}), by choosing $c_{5}$ (resp. $\alpha_{9}$) small enough (resp. large enough) in
(\ref{Eqn:InducStr3}).

\end{prop}


\vspace{3mm}
\subsection{The proof}
\label{Subsec:TheProofThm} $ $\\
We are now in position to prove Theorem \ref{Thm:Scat}. We define the following statement of induction
for $l \in \mathbb{N} $ \\

$\mathcal{P}(l)$: let
\begin{equation*}
\mathcal{C}_{l}  :=
\left\{
\begin{array}{ll}
w_{l},
\begin{array}{l}
w_{l} \quad \text{solution} \quad \text{of} \quad (\ref{Eqn:NlkgWdat}), \\
\exists \,  t_{l} \in \mathbb{R}^{+} \quad  s.t  \quad
E(Iw_{l}(t_{l})) \leq C_{in,2} A^{p+1} - 0.9 l c_{3} A^{-\alpha_{3}} ;
\end{array}
\end{array}
\right\},
\end{equation*}
then there exists $ \infty > L(l):=L_{N,s,A}(l)$ 
such that
\begin{equation}
\inf \{ \bar{C}, \, w_{l} \in \mathcal{C}_{l} \quad \text{and} \quad \| I w_{l} \|_{L_{t}^{2(p-1)-} L_{x}^{2(p-1)+} (\mathbb{R}^{+})} \leq \bar{C} \}  =L(l)
\label{Eqn:InducLtLx}
\end{equation}
and
\begin{equation}
\sup_{t \in \mathbb{R}^{+}} E(I w_{l}(t))  \leq ( C_{in,2} A^{p+1} - 0.9 l c_{3} A^{-\alpha_{3}} ) (1 + c_{6} A^{- \alpha_{10}}).
\label{Eqn:InducSupNrj}
\end{equation}
\\
We easily obtain that $\mathcal{P}(\bar{l} )$ holds for some $\bar{l} \lesssim A^{p+1 + \alpha_{3}}$ by applying Proposition
\ref{prop:EstSmallNrj}, choosing $N$ such that (\ref{Eqn:InducStr3}) holds. \\
\\
Our goal is then to show that if  $\mathcal{P}(l+1)$ holds, then $\mathcal{P}(l)$ also holds for $N$ and $s$ to be properly chosen. To this end  let $w_{l} \in \mathcal{C}_{l}$. Let $T>0$. Assume that (\ref{Eqn:InducLtLx}) restricted to $[0, t_{l} + T]$ holds for some $L(l) < \infty$ to be chosen. Choose $N$ such that (\ref{Eqn:InducStr2}), (\ref{Eqn:Cond4}), and (\ref{Eqn:InducStr3}) hold with $L_{w}$, $L_{w'}$ substituted respectively for
$L(l)$, $L(l+1)$. Clearly (\ref{Eqn:Cond4}) is the most constraining assumption to satisfy, choosing $C_{3}$ and $\alpha_{6}$ (resp. $c_{4}$)
large enough (resp. small enough). One may partition $[0, t_{l} + T]$ into subintervals $J$ such that $ \| I w_{l} \|_{L_{t}^{2(p-1)-} L_{x}^{2(p-1)+} (J)} =c_{2}$ (except maybe the last one), one may apply Proposition \ref{prop:LocalBd} and Proposition \ref{Prop:Acl} on each $J$, and then one may iterate to get
for $t \in [0,t_{l}+T]$

\begin{equation}
\begin{array}{ll}
\left| E(I w_{l}(t)) - E(I w_{l}(t_{l})) \right| & \lesssim
\frac{ ( C_{in,2} A^{p+1} - 0.9 l c_{3} A^{- \alpha_{3}})^{\frac{p+1}{2}} \langle L(l) \rangle^{2(p-1)-}}{N^{\frac{5-p}{2}-}} \\
& \leq  c_{6} A^{- \alpha_{10}} \left( C_{in,2} A^{p+1} - 0.9 l c_{3} A^{-\alpha_{3}} \right),
\end{array}
\nonumber
\end{equation}
using (\ref{Eqn:InducStr2}), and choosing $\alpha_{1}$ (resp. $c_{1}$) large (resp. small) enough. Therefore
(\ref{Eqn:InducSupNrj}) holds. It remains to prove (\ref{Eqn:InducLtLx}). To this end we let $M$ be such that

\begin{equation}
 C_{2} \frac{A^{\alpha_{4}}}{M^{\alpha_{5}}}  = \frac{1}{C_{4}^{ ( A \langle L(l+1) \rangle )^{\alpha_{7}}}},
\label{Eqn:DefM}
\end{equation}
Let $B >0$. If  $ \| I w_{l} \|_{L_{t}^{2(p-1)-} L_{x}^{2(p-1)+} ([0, t_{l} + T])} \geq 3B$, then we can find
$( \overline{t}_{l}, \overline{\overline{t}}_{l}) \in [0, t_{l}+ T]^{2}$
such that
$$
\| I w_{l} \|_{L_{t}^{2(p-1)-} L_{x}^{2(p-1)+} ([0, \overline{t}_{l}]) } \geq B, \quad
\| I w_{l} \|_{L_{t}^{2(p-1)-} L_{x}^{2(p-1)+} ([\overline{t}_{l}, \overline{\overline{t}}_{l}])}  \geq B, \quad \text{and} \quad
\| I w_{l} \|_{L_{t}^{2(p-1)-} L_{x}^{2(p-1)+} ([\overline{\overline{t}}_{l}, t_{l} + T])} \geq B.
$$
Assume that $ B \geq C_{1}^{( C_{1}^{A^{\alpha_{2}}}  M )^{ C_{1} A^{\alpha_{2}}}} $. Then, applying
Proposition \ref{Prop:SepLocMol} to  $J':=[\overline{t}_{l}, \overline{\overline{t}}_{l}]$, we see that there exists $T_{l} \in J^{'}$
and $w^{'}_{l}$ solution of (\ref{Eqn:NlkgWdat}) such that
\begin{align} \nonumber
E(Iw^{'}_{l}, T_{l}) & \leq ( C_{in,2} A^{p+1} - 0.9 l c_{3} A^{-\alpha_{3}}) (1+ c_{6} A^{-\alpha_{10}}) -
c_{3} A^{-\alpha_{3}} \\
& \leq C_{in,2} A^{p+1} -  0.9 (l+1) c_{3} A^{- \alpha_{3}},
\label{Eqn:WhyBeta}
\end{align}
choosing $c_{6}$ (resp. $\alpha_{10}$) small (resp. large) enough. Therefore $w^{'}_{l} \in \mathcal{C}_{l+1}$.
We deal with the case where (\ref{Eqn:FreeDecay}) holds with $J'':=[T_{l}, t_{l} + T]$ \footnote{ the other case $J'' := [t_{l},T_{l}] $ is
handled by a similar argument and therefore it is left to the reader.}. Applying
$\mathcal{P}(l+1)$ and Proposition \ref{prop:PerturbArg}, we see that
\begin{align*}
B \leq  \| I w_{l} \|_{L_{t}^{2(p-1)-} L_{x}^{2(p-1)+} ([\overline{\overline{t}}_{l}, t_{l} + T])}  \leq \| I w_{l} \|_{L_{t}^{2(p-1)-} L_{x}^{2(p-1)+} ([T_{l}, t_{l} + T])}
\lesssim  \left( \langle L(l+1) \rangle A \right)^{\alpha_{8}}.
\end{align*}
Therefore we see that $ \| I w_{l} \|_{L_{t}^{2(p-1)-} L_{x}^{2(p-1)+} ([0,t_{l}+T]) }  < \infty $, $\mathcal{P}(l)$ holds, and, moreover,
there exist $\alpha \gtrsim 1$ and $C \gtrsim 1$ such that

\begin{align}
L(l)   \lesssim  \max \left(  C_{1}^{ ( C_{1}^{A^{\alpha_{2}}} C_{2}^{\frac{1}{\alpha_{5}}} A^{\frac{\alpha_{4}}{\alpha_{5}}}
C_{4}^{\frac{ (A \langle L(l+1)  \rangle)^{\alpha_{7}}}{\alpha_{5}}} )^{ C_{1} A^{\alpha_{2}}} }, (\langle L(l+1) \rangle A)^{\alpha_{8}} \right)
\leq  C^{ C^{ \max^{\alpha}(A, \langle L(l+1) \rangle )^{C \max^{\alpha}(A, \langle L(l+1) \rangle) }}} .
\label{Eqn:LlInduc}
\end{align}
\\
\\
Iterating (\ref{Eqn:LlInduc}) $\bar{l}$ times, we see that for $0 \leq l \leq \bar{l}$ (even if it means increasing the value of $\alpha$
and $C$)

\begin{align*}
L(l) \leq \| L \|_{\infty}:= C^{C^{...^{ C^{\alpha}}}},
\end{align*}
where the height of the tower is $\sim A^{\alpha}$. Such an iteration is possible if $ \frac{ C_{3}^{ (A \langle \| L \|_{\infty} \rangle)^{\alpha_{6}}}
(A  \langle \| L \|_{\infty} \rangle)^{\alpha_{6}} } {N^{(1-s_{c})-}}  \leq c_{4}$. Pick $N$ such that
$ \frac{c_{4}}{2} \leq \frac{ C_{3}^{ (A \langle \| L \|_{\infty} \rangle)^{\alpha_{6}}}  (A \langle \| L \|_{\infty} \rangle) ^{\alpha_{6}} } {N^{(1-s_{c})-}}  \leq c_{4} $. From (\ref{Eqn:ChoiceNs}), we see that there exists $\tilde{s}$ such that $\mathcal{P}(l)$ holds for  $1  >  s > \tilde{s}$. Moreover $\tilde{s}$ can be chosen to be of the form (\ref{Eqn:Formtildes}) as $A \rightarrow \infty$.


\vspace{5mm}

\noindent \textbf{Global existence} \\
\\
From $\mathcal{P}(0)$ and (\ref{Eqn:BoundInitMolNrj}) we see that solutions of (\ref{Eqn:NlkgWdat}) with $1 > s > \tilde{s}$ and with data $(u_{0},u_{1}) \in H^{s} \times H^{s-1}$  satisfying
(\ref{Eqn:UpperBdInit}) satisfy

\begin{equation}
\sup_{t \in \mathbb{R^{+}}} E(Iu(t))   \lesssim A^{p+1}, \quad \text{and} \quad \| I u \|_{L_{t}^{2(p-1)-} L_{x}^{2(p-1)+} (\mathbb{R^{+}})}  < \infty.
\label{Eqn:BdLgEst}
\end{equation}
By time reversal symmetry, we may extend $\mathbb{R}^{+}$ to $\mathbb{R}$. By Plancherel, we have for all time $T \in \mathbb{R}$
$ \| ( u(T),\partial_{t} u(T)) \|^{2}_{H^{s} \times H^{s-1}}  \lesssim A^{p+1}$. This
proves global well-posedness. \footnote{Notice that global well-posedness was already proved in
\cite{bodaomiao} but, since it is a prerequisite to prove scattering, we reprove it. } \\
\\
\noindent \textbf{Global estimates} \\
\\
Now we claim that
\begin{equation}
Z_{m,s}(\mathbb{R},u)  \lesssim_{A} 1
\label{Eqn:GlobEst}
\end{equation}
for all $0 \leq m \leq s$. Indeed, by (\ref{Eqn:BdLgEst}), we may divide $[0, \infty)$ in subintervals $J:=[a,b]$ such that $\| Iu \|_{L_{t}^{2(p-1)-} L_{x}^{2(p-1)+}(J)} = c_{2}$ (except maybe the last one). Moreover, plugging $\langle D \rangle^{1-m} I $ into (\ref{Eqn:StrNlkg}), and by (\ref{Eqn:BdLgEst}), we have \\
\begin{equation}
\begin{array}{ll}
Z_{m,s}(J,u) &  \lesssim  \left\| ( \langle D \rangle I u(a), \partial_{t} I u(a) )  \right\|_{L^{2} \times L^{2}} +
\left\| \langle D \rangle^{1-m} I (|u|^{p-1} u) \right\|_{L_{t}^{\frac{2}{1+m}} L_{x}^{\frac{2}{2-m}} (J) }  \\
& \lesssim  A^{\frac{p+1}{2}} + \left\| \langle D \rangle^{1-m} I u \right\|_{L_{t}^{\frac{2}{m}+} L_{x}^{\frac{2}{1-m}-} (J)}
\left\| |u|^{p-1} \right\|_{L_{t}^{2-} L_{x}^{2+} (J)} \\
& \lesssim A^{\frac{p+1}{2}}  + Z_{m,s}(J,u)  \left(  \| P_{\ll N} u \|^{p-1}_{L_{t}^{2(p-1)-} L_{x}^{2(p-1)+} (J)} + \| P_{\gtrsim N} u \|^{p-1}_{L_{t}^{2(p-1)-} L_{x}^{2(p-1)+} (J)}  \right)  \\
& \lesssim   A^{\frac{p+1}{2}} + Z_{m,s}(J,u) \left( \| I u \|^{p-1}_{L_{t}^{2(p-1)-} L_{x}^{2(p-1)+} (J)} +
\frac{ \| \langle D \rangle^{1- s_{c}} I u \|^{p-1}_{L_{t}^{2(p-1)-} L_{x}^{2(p-1)+} (J)}}{N^{\frac{5-p}{2}-}}  \right)  \\
& \lesssim  A^{\frac{p+1}{2}} + Z_{m,s}(J,u) \left( c_{2}^{p-1} +
\frac{ Z^{p-1}_{s_{c},s} (u,J)} {N^{\frac{5-p}{2}-}}  \right) \\
\end{array}
\label{Eqn:Comput}
\end{equation}
By our choice of $N$, we have  $\frac{A^{\frac{p(p+1)}{2}}}{N^{\frac{5-p}{2}-}} \ll A^{\frac{p+1}{2}}$. Therefore a continuity
argument (first for $m=s_{c}$, then for the other $m$), we see that $Z_{m,s}(J,u) \lesssim A^{\frac{p+1}{2}} $. Iterating over $J$, we get (\ref{Eqn:GlobEst}).
\vspace{5 mm}

\noindent \textbf{Scattering} \\

Let $ \mathbf{v}(t)  := \left( u(t), \, \partial_{t} u(t) \right) $,
$ \mathbf{v}_{0} : = \left( u_{0}, \, u_{1} \right)$ and
\begin{equation*}
\mathbf{u_{nl}}(t)  := \left(
\begin{array}{l}
- \int_{0}^{t} \frac{\sin{\left((t-t^{'}) \langle D \rangle \right)}}{\langle D \rangle} \left( |u|^{p-1}(t^{'}) u(t^{'}) \right) \, d t^{'} \\ \\
- \int_{0}^{t} \cos{ \left( (t-t^{'}) \langle D \rangle \right)} \left( |u|^{p-1}(t^{'}) u(t^{'}) \right) \, d t^{'}
\end{array}
\right) \, .
\end{equation*}
Then we get from (\ref{Eqn:StrongSol}) $ \mathbf{v}(t)  = K(t) \mathbf{v_{0}} + \mathbf{u_{nl}}(t)$. Recall that the
solution $u$  scatters in $H^{s} \times H^{s-1}$  if there exists
$ \mathbf{v_{+,0}}  :=  (u_{+,0}, u_{+,1}) $ such that $ \left\|   \mathbf{v}(t)
 -K(t) \mathbf{v_{+,0}} \right\|_{H^{s} \times H^{s-1}} \rightarrow 0$  as $t \rightarrow \infty$. In other words, since
 $K$ is bounded on $H^{s} \times H^{s-1}$, it suffices to prove that the quantity
 $ \left\| K^{-1}(t) \mathbf{v}(t) -  \mathbf{v_{+,0}}   \right\|_{H^{s} \times H^{s-1}} \rightarrow 0$ as $t \rightarrow \infty$. A computation
shows that
\begin{equation*}
K^{-1}(t)  = \left(
\begin{array}{cc}
\cos{(t \langle D \rangle )} & -\frac{\sin{(t \langle D \rangle)}}{\langle D \rangle} \\ & \\
\langle D \rangle \sin {(t \langle D \rangle)} & \cos{(t \langle D \rangle)}
\end{array}
\right) \, .
\end{equation*}
But $ K^{-1}(t) \mathbf{v}(t)   =  \mathbf{v_{0}} - K^{-1}(t) \mathbf{u_{nl}}(t) $ and, by dualizing the Strichartz estimate
$\| e^{i t \langle D \rangle} f \|_{L_{t}^{\frac{2}{1-s}} L_{x}^{\frac{2}{1+s}}} \lesssim \| f \|_{H^{1-s}}$ (see
Proposition \ref{prop:StrEstLtqLxr}), we have
\begin{align*}  \nonumber
\| K^{-1}(t_{1}) \mathbf{u_{nl}}(t_{1}) - K^{-1}(t_{2}) \mathbf{u_{nl}}(t_{2}) \|_{H^{s} \times H^{s-1}}
& \lesssim \| |u|^{p-1} u \|_{L_{t}^{\frac{2}{1+s}} L_{x}^{\frac{2}{2-s}} ([t_{1}, \, t_{2}]) } \\
& \lesssim \| \langle D \rangle^{1-s} I \left(  |u|^{p-1} u  \right) \|_{L_{t}^{\frac{2}{1+s}} L_{x}^{\frac{2}{2-s}} ([t_{1}, \, t_{2}]) } \, .
\end{align*}
But, plugging $\langle D \rangle^{1-s} I$ into (\ref{Eqn:StrNlkg}) and modifying slightly (\ref{Eqn:Comput}), we get
\begin{align*}
\nonumber
\left\| \langle D \rangle^{1-s} I (|u|^{p-1} u) \right\|_{L_{t}^{\frac{2}{1+s}} L_{x}^{\frac{2}{2-s}} ([t_{1},t_{2}])} &
& \lesssim  Z_{s,s}([t_{1},t_{2}],u) \left( \| I u \|^{p-1}_{L_{t}^{2(p-1)-} L_{x}^{2(p-1)+} ([t_{1},t_{2}]) } + \frac{Z^{p-1}_{s_{c},s}
([t_{1},t_{2}],u)} {N^{\frac{5-p}{2}-}} \right) \nonumber .
\end{align*}
Therefore, from (\ref{Eqn:BdLgEst}) and (\ref{Eqn:GlobEst}), we see that the Cauchy criterion
is satisfied by $K^{-1}(t) v(t)$  and we conclude that $K^{-1}(t) v(t)$ has a limit in $H^{s} \times H^{s-1}$ as $t$ goes to
infinity. Moreover $ \lim \limits_{t \rightarrow \infty} \left\|  \mathbf{v}(t) -  K(t)\mathbf{v_{+,0}}   \right\|_{H^{s} \times H^{s-1}}
 = 0 $, with $\mathbf{v_{+,0}} := \left( u_{+,0}, u_{+,1} \right)$ given explicitly by
\begin{align*}
u_{+,0} & := u_{0} + \int_{0}^{\infty} \frac{\sin{(t^{'} \langle D \rangle)}}{\langle D \rangle} \left( |u|^{p-1}(t^{'}) u(t^{'}) \right) \, d
t^{'} ,\\
u_{+,1} & := u_{1} - \int_{0}^{\infty} \cos{(t^{'} \langle D \rangle)} \left( |u|^{p-1}(t^{'}) u(t^{'}) \right) \, d t^{'} \, .
\end{align*}

\section{Proof of local boundedness} 
\label{Sec:LocalBd}
\noindent In this section we prove Proposition \ref{prop:LocalBd}.
Plugging $\langle D \rangle ^{1-m}I$ into (\ref{Eqn:StrNlkg}), we have (with $\bar{J}':=[\bar{a}',\bar{b}']$)
\begin{equation}
Z_{m,s}(J^{'},w)  \lesssim  \| \langle D \rangle I w(\bar{a}') \|_{L^{2}} + \| \partial_{t} I w(\bar{a}') \|_{L^{2}} + \| \langle D \rangle^{1-m} I ( |w|^{p-1} w ) \|_{L_{t}^{\frac{2}{1+m}} L_{x}^{\frac{2}{2-m}}(\bar{J}^{'})},
\label{Eqn:PlugNrjStr}
\end{equation}
by (\ref{Eqn:InducNrjw}). There are three cases:

\begin{itemize}

\item $ m \leq s$. By slightly modifying (\ref{Eqn:Comput})

\begin{align*}
\nonumber Z_{m,s}(\bar{J}^{'},w)  & \lesssim  X^{\frac{1}{2}}_{w} + Z_{m,s}(\bar{J}^{'},w) \left( \| I w \|^{p-1}_{L_{t}^{2(p-1)-} L_{x}^{2(p-1)+} (\bar{J}^{'})} + \frac{ \| \langle D \rangle^{1- s_{c}} I w \|^{p-1}_{L_{t}^{2(p-1)-} L_{x}^{2(p-1)+} (\bar{J}^{'})}}{N^{\frac{5-p}{2}-}}  \right)  \\
& \lesssim  X^{\frac{1}{2}}_{w} + Z_{m,s}(\bar{J}^{'},w) \left( c_{2}^{p-1} + \frac{ Z^{p-1}_{s_{c},s}(\bar{J}^{'},w)}{N^{\frac{5-p}{2}-}}  \right)
\end{align*}
Again, choosing $\alpha_{1}$ (resp. $c_{1}$) large enough (resp. small enough) in (\ref{Eqn:InducStr2}), we have
$\frac{A^{\frac{p(p+1)}{2}}}{N^{\frac{5-p}{2}-}} << A^{\frac{p+1}{2}}$. Therefore a continuity  argument (first for $m=s_{c}$, then for the other $m$)
shows that $ Z_{m,s}(\bar{J}^{'},w) \lesssim  X^{\frac{1}{2}}_{w} $. \\

\item $m=1$. We estimate
\begin{align*}
\nonumber Z_{m,s}(\bar{J}^{'},w) & \lesssim   X^{\frac{1}{2}}_{w} + \|  I (|w|^{p-1} w) \|_{L_{t}^{1} L_{x}^{2} (\bar{J}^{'})}  \\
\nonumber & \lesssim  X^{\frac{1}{2}}_{w} + \| |w|^{p-1} w \|_{L_{t}^{1} L_{x}^{2}(\bar{J}^{'})}  \\
 & \lesssim X^{\frac{1}{2}}_{w} + \| |P_{\ll N} w|^{p-1}  P_{\ll N} w  \|_{L_{t}^{1} L_{x}^{2} (\bar{J}^{'})} +  \| |P_{\ll N} w|^{p-1}
P_{\gtrsim N} w \|_{L_{t}^{1} L_{x}^{2}(\bar{J}^{'})}  \\
& + \| |P_{\gtrsim N} w|^{p-1} P_{<<N} w \|_{L_{t}^{1} L_{x}^{2} (\bar{J}^{'})} +  \| |P_{\gtrsim N} w|^{p-1} P_{\gtrsim N} w \|_{L_{t}^{1}
L_{x}^{2}(\bar{J}^{'})}  \\
\nonumber & \lesssim  X^{\frac{1}{2}}_{w} + B_{1} + B_{2} + B_{3} + B_{4}  \, .
\end{align*}
We estimate
\begin{align*}
B_{1} & \lesssim  \| I w \|^{p-1}_{L_{t}^{2(p-1)-} L_{x}^{2(p-1)+} (\bar{J}^{'})} \| I w \|_{L_{t}^{2+} L_{x}^{\infty-}(\bar{J}^{'})}  \\
& \lesssim  c_{2}^{p-1}  Z_{1,s}(w,\bar{J}^{'}),
\end{align*}

\begin{align*}
B_{2} & \lesssim \| I w \|^{p-1}_{L_{t}^{p-1} L_{x}^{\frac{6(p-1)}{p-3}} (\bar{J}^{'}) } \| P_{\gtrsim N} w \|_{L_{t}^{\infty}
L_{x}^{\frac{6}{6-p}}(\bar{J}^{'})}  \\
& \lesssim  Z^{p-1}_{1,s}(\bar{J}^{'},w) \frac{\| \langle D \rangle I w \|_{L_{t}^{\infty} L_{x}^{2}(\bar{J}^{'})} }{N^{\frac{5-p}{2}-}} \\
& \lesssim \frac{Z^{p-1}_{1,s}(\bar{J}^{'},w)}{N^{\frac{5-p}{2}-}} X_{w}^{\frac{1}{2}},
\end{align*}
\begin{align*}
B_{3} & \lesssim  \frac{ \| \langle D \rangle^{1-s_{c}} I w \|^{p-1}_{L_{t}^{2(p-1)-} L_{x}^{2(p-1)+}(\bar{J}^{'})}}{N^{\frac{5-p}{2}-}} \| Iw
\|_{L_{t}^{2+} L_{x}^{\infty -}(\bar{J}^{'})} \\
& \lesssim \frac{Z^{p-1}_{s_{c},s}(\bar{J}^{'},w)}{N^{\frac{5-p}{2}-}} Z_{1,s}(\bar{J}^{'},w)
\lesssim  \frac{X_{w}^{\frac{p-1}{2}} }{N^{\frac{5-p}{2}-}} Z_{1,s}(\bar{J}^{'},w),
\end{align*}

\begin{align*}
B_{4} &  \lesssim \| P_{\gtrsim N} w \|^{p}_{L_{t}^{p} L_{x}^{2p}(\bar{J}^{'})} \lesssim  \frac{ \| \langle D \rangle^{1- \left( \frac{3p-5}{2p} \right)} I w \|^{p}_{L_{t}^{p} L_{x}^{2p}(\bar{J}^{'})}}{N^{\frac{5-p}{2}-}} \\
& \lesssim  \frac{ Z^{p}_{\frac{3p-5}{2p},s} (\bar{J}^{'},w) }{N^{\frac{5-p}{2}-}} \lesssim \frac{X_{w}^{\frac{p}{2}}}{N^{\frac{5-p}{2}-}} \, .
\end{align*}
Therefore, since again $N$ satisfies (\ref{Eqn:InducStr2}), we see by  a continuity argument that $ Z_{1,s}(\bar{J}^{'},w) \lesssim X^{\frac{1}{2}}_{w} $.

\item $s<  m < 1$: $ Z_{m,s}(J^{'},w) \lesssim X^{\frac{1}{2}}_{w} $ follows by interpolating between $m=s$ and $m=1$.

\end{itemize}

\section{Proof of separation of the localized mollified energy} 
\label{Sec:SepLocMol}
\noindent In this section we prove Proposition \ref{Prop:SepLocMol}. The proof of Proposition \ref{Prop:SepLocMol} relies upon three lemmas
that we show in the next subsections. 

\subsection{Lemmas \ref{lem:ConcNrj}, \ref{lem:UseMorEst}, and \ref{lem:StrVsConcNrj}}

The first lemma shows that if there is concentration of the target norm of the solution on a subinterval of $J^{'}$
in the sense of (\ref{Eqn:ConcPot}), then this also means that the potential term of the mollified energy and the size of this subinterval are substantial.

\begin{lem}
Assume that
\begin{equation}
\| I w \|_{L_{t}^{2(p-1)-} L_{x}^{2(p-1)+} (\bar{J}^{'}) }  = c_{2}.
\label{Eqn:ConcPot}
\end{equation}
Then there exist a subinterval $\bar{K}^{'} \subset \bar{J}^{'}$, a
number $R \geq 1$, a point $\bar{x}^{'} \in  \mathbb{R}^{3}$ and constants $C_{5}$, $c_{7}$, $c_{8}$, $\alpha_{11}$,...,
$\alpha_{13}$ such that



\begin{equation}
\begin{array}{ll}
R:= C_{5} A^{\alpha_{11}}, \, |\bar{K}^{'}| & = c_{7} A^{-\alpha_{12}},
\end{array}
\label{Eqn:SizeK}
\end{equation}
and for all $t \in \bar{K}^{'}$
\begin{equation}
\int_{|x-\bar{x}^{'}| \leq R} |Iw (t,x)|^{p+1} \, dx  \geq  c_{8}  A^{-\alpha_{13}}.
\label{Eqn:ConcNrjMol}
\end{equation}


\label{lem:ConcNrj}
\end{lem}
The second lemma shows that if we consider a partition of $J'$ into subintervals where the
target norm of the solution concentrates, then these subintervals must be large on average.
In order to prove this lemma, we shall mostly use the previous lemma and the Almost Morawetz-Strauss estimate
(\ref{Eqn:AlmMorStrEst}).

\begin{lem}

Let $ (\bar{J}^{'}_{j}=[\bar{a}^{'}_{j},\bar{a}^{'}_{j+1}])_{ 1 \leq j \leq \bar{j}} $ be a partition of $J^{'}$ such that
$ \| Iw  \|_{L_{t}^{2(p-1)-} L_{x}^{2(p-1)+} (\bar{J}^{'}_{j})} = c_{2} $, except maybe the last one. Then there exist $\bar{t}^{'}_{j} \in \bar{J}^{'}_{j}$
and a constant $\alpha_{14}$ such that
\begin{equation}
\sum_{j=1}^{\bar{j}-1} \frac{1}{\bar{t}^{'}_{j} + 1}  \lesssim A^{\alpha_{14}}.
\label{Eqn:RepartJj}
\end{equation}

\label{lem:UseMorEst}
\end{lem}
The third lemma shows that if the target norm of the solution is too large in the sense of (\ref{Eqn:UpperBdIu}) then we can find a large subinterval where some norms are small compare with the concentration of mollified energy in the sense of (\ref{Eqn:BoundZpr})

\begin{lem}

Let $M \geq 1$. Then, there exist $C_{6}$, $c_{9}$, $\alpha_{15}$, and $\alpha_{16}$  for all $ 0< \epsilon \leq 1$, there exist $R^{'} \in (1,\infty)$, $\bar{x} \in \mathbb{R}^{3}$ and $J'' := [S,T] $ (or $J'' :=[T,S]$) such that $ J'' \subset J^{'}$ and if


\begin{equation}
\| I w \|_{L_{t}^{2(p-1)-} L_{x}^{2(p-1)+} (J^{'})} \geq
C_{6}^{ (C_{6}^{( A \epsilon^{-1})^{\alpha_{15}}} M )^{ C_{6} A^{\alpha_{15}}}},
\label{Eqn:UpperBdIu}
\end{equation}
then
\begin{equation}
\tilde{Z}(J'',w) \leq c_{9} A^{-\alpha_{16}} \leq E( Iw(S),B(\bar{x},R^{'})),
\label{Eqn:BoundZpr}
\end{equation}

\begin{equation}
|J''| \geq M R^{'},
\label{Eqn:LowerBoundST}
\end{equation}
and

\begin{equation}
\left\| \frac{Iw(S)}{\langle x-\bar{x} \rangle} \right\|_{L^{2}} \leq \epsilon \, .
\label{Eqn:UpBdWeight}
\end{equation}

\label{lem:StrVsConcNrj}
\end{lem}

\subsection{The proof}

\noindent We may assume without loss of generality that $S < T$. We apply Lemma \ref{lem:StrVsConcNrj} with
$\epsilon  <<  \min \left( A^{-\frac{\alpha_{16}}{2}}, A^{- \alpha_{16} - \frac{p+1}{2}} \right) $. Notice that with this choice of
$\epsilon$, the condition (\ref{Eqn:UpperBdIu}) becomes (\ref{Eqn:TooLarge}), choosing $\alpha_{2}$ large enough. The proof is made of several steps: \\
\\
Step 1. \textit{Construction of the free Klein-Gordon equation $v$ and proof of (\ref{Eqn:ConservkinetIv})}. \\
\\
Let $P(y):= \left\{ y \in \mathbb{R}^{3}, \, E \left( Iw(S), B (y,1) \right)  \leq  c_{9} A^{-\alpha_{16}} \right\}$.
Then by (\ref{Eqn:InducNrjw}) there exists $\bar{\bar{x}} \in \mathbb{R}^{3} $ such that
$ |\bar{x} - \bar{\bar{x}}| \lesssim A^{p+ 1 + \alpha_{16}} $ and $ P(\bar{\bar{x}}) $ is  true. Hence, using
also (\ref{Eqn:BoundZpr}) we see that there exists a constant $C$ and $\Gamma \in [\frac{1}{2}, R^{'} +  C A^{p+1+ \alpha_{16}} ]$ such
that

\begin{equation}
E \left( I w(S), B(\bar{\bar{x}}, \Gamma)  \right)  = c_{9} A^{-\alpha_{16}} \, .
\end{equation}
Let $v$ be the solution of the free Klein-Gordon equation with data

\begin{equation*}
\begin{cases}
v(S)  := I^{-1} \left( \chi \left( \frac{x - \bar{\bar{x}}}{\Gamma} \right) Iw(S) \right) \\
\partial_{t} v (S)  := I^{-1} \left( \chi \left( \frac{x - \bar{\bar{x}}}{\Gamma} \right)   \partial_{t} I w(S) \right),
\end{cases}
\end{equation*}
where $\chi$ is a smooth function such that $\chi(x)=1$ if $|x| \leq 1$ and $\chi(x)=0$ if $|x| \geq 2$. By (\ref{Eqn:InducNrjw}) and
(\ref{Eqn:UpBdWeight}) we see that there exists a constant $C$ such that

\begin{align}
E_{c} \left( Iv(S) \right) & \leq E \left( Iw(S), B(\bar{\bar{x}}, \Gamma ) \right) +  \frac{C}{\Gamma^{2}} \int_{
\Gamma \leq |x-
\bar{\bar{x}}| \leq 2 \Gamma} |Iw(S)|^{2} \, dx  + \nonumber \\ \label{Eqn:IneqNrjS}
&  \frac{C}{\Gamma} \int_{ \Gamma \leq |x - \bar{\bar{x}}| \leq 2 \Gamma}
|Iw(S)| |\nabla Iw (S)| \, dx  \\ \nonumber
& \leq  c_{9} A^{-\alpha_{16}} + C  A^{\frac{p+1}{2}} \left\| \frac{I w(S)}{\langle x - \bar{\bar{x}} \rangle}
\right\|_{L^{2}}
 + C  \left\| \frac{I w(S)}{\langle x - \bar{\bar{x}} \rangle} \right\|^{2}_{L^{2}} \\ \nonumber
& \lesssim A^{-\alpha_{16}},
\end{align}
Hence, using also the conservation of $E_{c}(Iv)$, we see that (\ref{Eqn:ConservkinetIv}) holds. \\
\\
Step 2. \textit{Proof of the decay (\ref{Eqn:FreeDecay})}. \\
\\
\label{Sec:FreeDecay}
\vspace{2mm}
By interpolation we see that one can choose one can choose $m < s_{c}$ close to $s_{c}$ and $ \gamma > 0$ close to zero such that

\begin{align}
\| Iv \|_{L_{t}^{2(p-1)-} L_{x}^{2(p-1)+} (T,b)} & \lesssim \| \langle D \rangle^{1-m} I v
\|^{1- \gamma}_{L_{t}^{\frac{8}{3-2m}-} L_{x}^{\frac{8}{3-2m}+} (T,b)} \| \langle  D \rangle^{- \left( \frac{3}{2}++ \right) } I
v \|^{\gamma}_{L_{t}^{\infty}
L_{x}^{\infty} (T,b)} \nonumber \\
& \lesssim  A^{- (1-\gamma)\alpha_{16}} \| I v \|^{\gamma}_{L_{t}^{\infty} B_{\infty,2}^{- \left( \frac{3}{2} + \right) } (T,b) } \nonumber \\
& \lesssim A^{- (1-\gamma)\alpha_{16}} \frac{1}{|T-S|^{\frac{3 \gamma}{2}}} \left( \| I v(S) \|^{\gamma}_{B^{1-}_{1,2}} +
\| \partial_{t} I v(S) \|^{\gamma}_{B^{0-}_{1,2}}
\right) \nonumber \\
& \lesssim  A^{- (1-\gamma)\alpha_{16}} \frac{1}{|T-S|^{\frac{3 \gamma}{2}}} \left( \Gamma^{\frac{3\gamma}{2}} \| I w(S) \|^{\gamma}_{L^{2}} +
\Gamma^{\frac{\gamma}{2}} \| \nabla I w (S) \|^{\gamma}_{L^{2}} + \Gamma^{\frac{3 \gamma}{2}} \| \partial_{t} I w (S) \|^{\gamma}_{L^{2}} \right) \nonumber \\
& \leq C_{2}  \frac{A^{\alpha_{4}}}{M^{\alpha_{5}}},
\label{Eqn:FreeDecayComput}
\end{align}
using also (\ref{Eqn:LowerBoundST}), (\ref{Eqn:IneqNrjS}) combined with (\ref{Eqn:StrNlkg}), and the following  dispersive estimate (see \cite{ginebvelo}, Lemma 2.1)
\begin{equation}
\|  e^{it \langle D \rangle } \phi \|_{B^{-\frac{5}{4}}_{\infty,2}}  \lesssim \frac{1}{|t|^{\frac{3}{2}}} \| \phi \|_{B^{\frac{5}{4}}_{1,2} }.
\label{Eqn:DispFreeKg}
\end{equation}
Step 3. \textit{Proof of the separation of the localized mollified energy (\ref{Eqn:SepMolNrj}).}\\
\\
Let $\bar{w}:=w-v$. Then

\begin{align}
E \left( I \bar{w}(S) \right)  & \leq \frac{1}{2} \int_{\mathbb{R}^{3}}  \left(  1 - \chi^{2} \left( \frac{x- \bar{\bar{x}}}{\Gamma}
\right) \right) |\nabla I w(S)|^{2} \, dx \nonumber \\
&  + \frac{1}{p+1} \int_{\mathbb{R}^{3}} \left( 1 - \chi^{p+1} \left( \frac{x-
\bar{\bar{x}}}{\Gamma} \right)
\right) |Iw(S)|^{p+1} \, dx  \nonumber \\ 
 &+ \frac{1}{2} \int_{\mathbb{R}^{3}}  \left(  1 - \chi^{2} \left( \frac{x- \bar{\bar{x}}}{\Gamma}  \right) \right) |\partial_{t} Iw(S)|^{2}
\,
dx + \frac{1}{2} \int_{\mathbb{R}^{3}} \left(  1 - \chi^{2} \left( \frac{x- \bar{\bar{x}}}{\Gamma}  \right) \right) |Iw(S)|^{2} \, dx \nonumber \\
 &+ \frac{C}{\Gamma} \int_{ \Gamma \leq |x-\bar{\bar{x}}| \leq 2 \Gamma}   | Iw(S)| | \nabla I w(S)| \, dx   + \frac{C}{\Gamma^{2}}
\int_{ \Gamma \leq |x - \bar{\bar{x}}| \leq 2 \Gamma} |Iw(S)|^{2} \, dx \nonumber  \\
 &\leq  \sup_{t \in J} E(Iw(t))  - c_{9}  A^{-\alpha_{16}} + C
 A^{\frac{p+1}{2}} \left\| \frac{I w(S)}{\langle x - \bar{\bar{x}}
\rangle} \right\|_{L^{2}}   + C  \left\| \frac{I w(S)}{\langle x - \bar{\bar{x}} \rangle} \right\|^{2}_{L^{2}} \nonumber \\
 &\leq  \sup_{t \in J} E(Iw(t)) - \frac{c_{9} A^{-\alpha_{16}}}{2} \, .
\label{IneqChainNrj}
\end{align}
Let

\begin{equation*}
\bar{Z}([S,T], f ) := \| I f \|_{L_{t}^{2(p-1)-} L_{x}^{2(p-1)+} ([S,T])} +
\| \partial_{t} \langle D \rangle^{-\frac{1}{2}} I f \|_{L_{t}^{4} L_{x}^{4} ([S,T])} +
\|  \langle D  \rangle^{\frac{1}{2}} I f \|_{L_{t}^{4+} L_{x}^{4-} ([S,T])}
\end{equation*}
Plugging $\langle D \rangle^{1- \frac{1}{2}} I$ into (\ref{Eqn:StrNlkg}), we see, by (\ref{Eqn:BoundZpr}) and the Sobolev
embedding

\begin{equation*}
\| I f \|_{L_{t}^{2(p-1)-} L_{x}^{2(p-1)+} ([S,T])} \lesssim
\| \langle D  \rangle^{1 - \frac{1}{2}}  I f \|_{L_{t}^{2(p-1)-} L_{x}^{\frac{6(p-1)}{2p-3}+} ([S,T]) },
\end{equation*}
that

\begin{align}
\label{Eqn:Zbarw}
 &\bar{Z}([S,T],\bar{w}) \lesssim \left\| \langle D \rangle^{1-\frac{1}{2}} I (|w|^{p-1} w)  \right\|_{L_{t}^{\frac{4}{3}} L_{x}^{\frac{4}{3}} ([S,T])} \\ \nonumber
& \lesssim \| \langle D \rangle^{1-\frac{1}{2}}  I w \|_{ L_{t}^{4+} L_{x}^{4-} ([S,T])}
\left( \|  P_{<<N} w  \|^{p-1}_{L_{t}^{2(p-1)-} L_{x}^{2(p-1)+} ([S,T])} + \|  P_{\gtrsim N} w  \|^{p-1}_{L_{t}^{2(p-1)-} L_{x}^{2(p-1)+} ([S,T])} \right) \\ \nonumber
& \lesssim  \| \langle D \rangle^{1-\frac{1}{2}}  I w \|_{ L_{t}^{4+} L_{x}^{4-} ([S,T])}
\left( \| I w \|^{p-1}_{L_{t}^{2(p-1)-} L_{x}^{2(p-1)+} ([S,T]) } +
\frac{ \| \langle  D \rangle^{1-\frac{1}{2}} I w \|^{p-1}_{L_{t}^{2(p-1)-} L_{x}^{\frac{6(p-1)}{2p-3}+}([S,T])  } }
{N^{\frac{5-p}{2}-}} \right)  \\ \nonumber
&\lesssim \tilde{Z}^{p}([S,T],w)  \\ \nonumber
& \leq \frac{c_{9} A^{-\alpha_{16}}}{1000}.
\end{align}
We compute
\begin{align*}
\partial_{t} E(I \bar{w}) & = \int_{\mathbb{R}^{3}} \Re \left( \overline{\partial_{t} I \bar{w}}
\left( \partial_{tt} I \bar{w}  - \triangle I \bar{w} + I \bar{w} + F(I \bar{w})  \right) \right) \, dx \\ \nonumber
& = \int_{\mathbb{R}^{3}} \Re \left( \overline{\partial_{t} I \bar{w}} ( F(I \bar{w}) - I F(w) ) \right) \, dx  \,.
\end{align*}
Now, we decompose $ E(I \bar{w}(T)) - E(I\bar{w}(S))  = X_{1,1} + X_{1,2} + X_{2} $ where
\begin{align}
 X_{1,1} &  : = \int_{S}^{T} \int_{\mathbb{R}^{3}} \Re{ \left( \overline{\partial_{t} I \bar{w}}  \left( F(w) -IF(w)  \right)
\right)} \,dx \, dt \, \nonumber \\
X_{1,2} &  := \int_{S}^{T} \int_{\mathbb{R}^{3}} \Re{ \left( \overline{\partial_{t} I \bar{w}}  \left( F(Iw) -F(w)  \right)
\right)} \, dx \, dt \nonumber \\
X_{2} & := \int_{S}^{T} \int_{\mathbb{R}^{3}} \Re{ \left( \overline{\partial_{t} I \bar{w}}  \left( F (I \bar{w}) - F(Iw)
\right) \right)} \, dx \, dt.
\label{Eqn:DefX1112}
\end{align}
We estimate $X_{2}$. From (\ref{Eqn:BoundZpr}) and (\ref{Eqn:Zbarw})
\begin{align*}
|X_{2}| & \lesssim \| \partial_{t} \langle D \rangle^{-\frac{1}{2}} I \bar{w}  \|_{L_{t}^{4} L_{x}^{4} ([S,T])} \| \langle D
\rangle^{\frac{1}{2}}
(F(I \bar{w}) -F(Iw) ) \|_{L_{t}^{\frac{4}{3}} L_{x}^{\frac{4}{3}} ([S,T])  } \\
& \lesssim \| \partial_{t} \langle D \rangle^{-\frac{1}{2}} I \bar{w}  \|_{L_{t}^{4} L_{x}^{4} ([S,T])}  \left( \| I w
\|^{p-1}_{L_{t}^{2(p-1)-}
L_{x}^{2(p-1)+} ([S,T])} + \| I \bar{w} \|^{p-1}_{ L_{t}^{2(p-1)-} L_{x}^{2(p-1)+} ([S,T]) } \right) \\
& \left( \| \langle D \rangle^{ 1 -\frac{1}{2}} I w \|_{L_{t}^{4+} L_{x}^{4-}([S,T])} + \| \langle D \rangle^{\frac{1}{2}} I
\bar{w} \|_{L_{t}^{4+} L_{x}^{4-}([S,T])}    \right) \\
& \leq \frac{c_{9} A^{- \alpha_{16}}}{1000}.
\end{align*}
Hence, using also (\ref{IneqChainNrj}) and Result \ref{Eqn:BoundX11} (see Appendix A) we see that (\ref{Eqn:SepMolNrj}) holds, with
$w'$ solution of (\ref{Eqn:NlkgWdat}) such that $w'(T):=\bar{w}(T)$.

\subsection{Proof of Lemma \ref{lem:ConcNrj}} $ $\\
The proof is made of four steps: \\
\\
Step 1. \textit{Lower bound of the size of $\bar{J}^{'}$.} \\
\\
We see from Proposition \ref{prop:LocalBd} that if $p \geq 4$, then

\begin{align*}
c_{2} & = \| I w \|_{L_{t}^{2(p-1)-} L_{x}^{2(p-1)+} (\bar{J}^{'})} \\
& \lesssim |\bar{J}^{'}|^{\frac{5-p}{2(p-1)} +}  \| \langle D \rangle^{1- \frac{1}{2}} I w \|_{L_{t}^{\frac{2(p-1)}{p-4}}
L_{x}^{\frac{6(p-1)}{p+2}} (\bar{J}^{'})} \\
& \lesssim |\bar{J}^{'}|^{\frac{5-p}{2(p-1)} +} A^{\frac{p+1}{2}},
\end{align*}
and if $p < 4$, then
\begin{align*}
c_{2} & = \| I w \|_{L_{t}^{2(p-1)-} L_{x}^{2(p-1)+} (\bar{J}^{'})} \\
& \lesssim  |\bar{J}^{'}|^{\frac{1}{2(p-1)}+ }  \| \langle D \rangle^{1- \frac{1}{2} } I w  \|_{L_{t}^{\infty} L_{x}^{3} (\bar{J}^{'})} \\
& \lesssim |\bar{J}^{'}|^{\frac{1}{2(p-1)}+} A^{\frac{p+1}{2}} \, .
\end{align*}
Therefore, we conclude that there exists a constant $c$  such that
\begin{equation}
|\bar{J}^{'}| \geq c \times
\left\{
\begin{array}{l}
A ^{-\frac{(p+1)(p-1)}{5-p}+}, \, p \geq 4 \\
A^{-(p+1)(p-1)+}, \, p < 4
\end{array}
\right.
\label{Eqn:LowerBoundJ}
\end{equation}\\
Step 2. \textit{Lower bound of $\| P_{M} w \|_{L_{t}^{\infty} L_{x}^{\infty}(\bar{J}^{'})}$ for some $M \in 2^{\mathbb{N}}$ }  \\
\\
From Proposition \ref{prop:LocalBd}

\begin{align*}
c_{2} = & \| I w \|_{L_{t}^{2(p-1)-} L_{x}^{2(p-1)+}(\bar{J}^{'})} \\
&  \lesssim \| \langle D \rangle^{1- \left (\frac{1}{2} + \right)} I w \|^{\frac{2}{p-1}}_{L_{t}^{4-}
L_{x}^{4+} (\bar{J}^{'})} \| \langle D \rangle^{ \left( - \frac{1}{p-3} \right) +}  I w \|^{1- \frac{2}{p-1}}_{L_{t}^{\infty} L_{x}^{\infty} (\bar{J}^{'})} \\
& \lesssim A^{\frac{p+1}{p-1}} \| \langle D \rangle^{ \left( - \frac{1}{p-3} \right) +}  I w \|^{1-
\frac{2}{p-1}}_{L_{t}^{\infty} L_{x}^{\infty} (\bar{J}^{'})}.
\end{align*}
Thus we have $ \| \langle D \rangle^{ \left( - \frac{1}{p-3}  \right) +}  I w  \|_{L_{t}^{\infty} L_{x}^{\infty} (\bar{J}^{'})}
\gtrsim A^{-\frac{p+1}{p-3}} $  and, by the pigeonhole principle, we conclude that there exists $M \in 2^{\mathbb{N}}$ such that
\begin{equation}
\| P_{M} I w \|_{L_{t}^{\infty} L_{x}^{\infty}(\bar{J}^{'})} \gtrsim \langle M \rangle^{\frac{1}{p-3}-}  A^{-\frac{p+1}{p-3}} \,.
\label{Eqn:UpBdPM1}
\end{equation}
On the other hand,
\begin{equation}
\| P_{M} I w \|_{L_{t}^{\infty} L_{x}^{\infty} (\bar{J}^{'})}  \lesssim   \langle M \rangle^{\frac{1}{2}} \| \langle D \rangle I w \|_{L_{t}^{\infty}
L_{x}^{2} (\bar{J}^{'})}
 \lesssim \langle M \rangle^{\frac{1}{2}}  A^{\frac{p+1}{2}} .
\label{Eqn:UpBdPm2}
\end{equation}
Combining (\ref{Eqn:UpBdPM1}) and (\ref{Eqn:UpBdPm2}), we see that
\begin{equation}
\langle M \rangle \lesssim  A^{\frac{(p+1)(p-1)}{5-p}+} \, .
\label{Eqn:UpperBdM}
\end{equation}
\\
Step 3. \textit{Control of $| P_{M}  I w  (t,\bar{x}^{'}) |$ for some $\bar{x}^{'} \in \mathbb{R}^{3}$ and for all
$t \in \bar{K}^{'}$, $\bar{K}^{'} \subset \bar{J}^{'}$ to be defined
shortly.}  \\
\\
By (\ref{Eqn:UpBdPM1}), there exists $(\bar{t}^{'},\bar{x}^{'})$ such that
\begin{equation}
|P_{M} I w(\bar{t}^{'},\bar{x}^{'})| \gtrsim \langle M \rangle^{\frac{1}{p-3}-}  A^{-\frac{p+1}{p-1}} \, .
\label{Eqn:aux1}
\end{equation}
But, by (\ref{Eqn:InducNrjw}) and (\ref{Eqn:UpperBdM}), we see that
\begin{align}
| P_{M} I w(t,\bar{x}^{'}) - P_{M} I w(\bar{t}^{'},\bar{x}^{'})| & \lesssim \sup_{s \in (\bar{t}{'},t)} \| \partial_{s} I w(s) \|_{L_{x}^{2}} \langle M
\rangle^{\frac{3}{2}}
|t - \bar{t}^{'}| \nonumber  \\
& \lesssim  A^{\frac{(p+1)^{2}}{5-p}}  |t - \bar{t}^{'}| \, ,
\label{Eqn:aux2} \\ \nonumber
\end{align}
Therefore, in view of (\ref{Eqn:LowerBoundJ}), (\ref{Eqn:aux1}) and (\ref{Eqn:aux2}), choosing $c_{7}$ (resp. $\alpha_{12}$) small enough
(resp. large enough), we see that either $ [\bar{t}^{'}, \bar{t}^{'} + c_{7} A^{-\alpha_{12}}] \subset \bar{J}^{'}$
(in this case we let $\bar{K}^{'}:=[\bar{t}^{'}, \bar{t}^{'} + c_{7} A^{-\alpha_{12}}] $), or
$ [\bar{t}^{'} - c_{7} A^{-\alpha_{12}} , \bar{t}^{'}] \subset \bar{J}^{'}$ (in this case
we let $\bar{K}^{'}:=[\bar{t}^{'} - c_{7} A^{-\alpha_{12}}, \bar{t}^{'}] $ ) and
\begin{equation}
|P_{M} Iw(t,\bar{x}^{'})|  \gtrsim  A^{-\frac{p+1}{p-1}}, \, t \in \bar{K}^{'}
\label{Eqn:PM}
\end{equation}
\\
Step 4. \textit{Lower bound of potential mollified energy.} \\
\\
Let $R>0$  to be fixed shortly. Let $\Psi:= \psi$ if $M > 1$ and
$\Psi:= \phi$ is $M=1$. By (\ref{Eqn:PM}) we have
\begin{equation*}
M^{3}(B_{1}+B_{2})  \gtrsim  A^{-\frac{p+1}{p-1}}
\end{equation*}
where $ B_{1}:= \int_{|y| \leq R} |\check{\Psi}(My)| |Iw(t,\bar{x}^{'}-y)| \, dy $ and
$ B_{2}:=  \int_{|y| \geq R} |\check{\Psi}(My)| |Iw(t,\bar{x}^{'}-y)| \, dy $. We have
\begin{align*}
B_{1}  & \lesssim \left( \int_{|y| \leq R} |\check{\Psi} (My) |^{\frac{p+1}{p}} \, dy \right)^{\frac{p}{p+1}} \left( \int_{|y| \leq R}
|Iw(t,\bar{x}^{'}-y)|^{p+1}
\, dy \right)^{\frac{1}{p+1}} \\
& \lesssim  \left( \int_{|y-\bar{x}^{'}| \leq R} |Iw(t,y)|^{p+1} \, dy \right)^{\frac{1}{p+1}} M^{-\frac{3p}{p+1}}
\end{align*}
and
\begin{align*}
B_{2} & \lesssim M^{-\frac{3}{2}} \| \check{\Psi} \|_{L^{2}(|y| \geq M R)} \| I w (t) \|_{L_{x}^{2}} \\
& \lesssim  M^{-\frac{3}{2}} \| \check{\Psi} \|_{L^{2}(|y| \geq M R)}   A^{\frac{p+1}{2}} \, .
\end{align*}
The fast decay of $\check{\Psi}$ implies  $$ \| \check{\Psi} \|_{L^{2}(|y| \geq  MR)} \lesssim \frac{1}{(MR)^{\frac{3}{2}}}.$$ Hence, if
$ R  := C_{5} A^{\alpha_{11}}$ (with $C_{5}$ and $\alpha_{11}$ large enough), then
\begin{equation*}
\int_{|y-\bar{x}^{'}| \leq R} |Iw(t,y)|^{p+1} \, dy  \geq c_{8}  A^{-\alpha_{13}},
\end{equation*}
for all $t \in \bar{K}^{'}$.

\subsection{Proof of Lemma \ref{lem:UseMorEst}} $ $\\
Let $j \in [1,\bar{j}-1]$. Recall that, by Lemma \ref{lem:ConcNrj}, on each $\bar{J}^{'}_{j}$, there exist $\bar{x}^{'}_{j} \in \mathbb{R}^{3}$ and
$\bar{K}^{'}_{j}=[\bar{t}^{'}_{j},\bar{t}^{'}_{j+1}] \subset \bar{J}^{'}_{j}$ such
that
\begin{equation}
\int_{|x- \bar{x}^{'}_{j}| \leq R } |I w(t,x)|^{p+1} \, dx  \geq c_{8} A^{-\alpha_{13}}
\label{Eqn:LowerPotNrj}
\end{equation}
for all $t \in \bar{K}^{'}_{j}$, with $R= C_{5} A^{\alpha_{11}}$  and
\begin{equation}
|\bar{K}^{'}_{j}|  = c_{7} A^{-\alpha_{12}}.
\label{Eqn:SizeKj}
\end{equation}


 We construct (see \cite{nakanash1}) a set
$\mathcal{P}:= \{ j_{1},....,j_{l} \} \subset \{ 1,...,\bar{j} -1 \}$. Initially $j_{1}=1$. Then let $j_{k+1}$ be the minimal $j$ such that
\begin{equation}
|\bar{x}^{'}_{j} - \bar{x}^{'}_{j_{k+1}}|  \geq |\bar{t}^{'}_{j}- \bar{t}^{'}_{j_{k+1}}| + 100R
\label{Eqn:Boundxjxk}
\end{equation}
for $j:=j_{1},...,j_{k}$. Observe that $J^{'} = \bigcup_{j_{k} \in \mathcal{P}} A^{'}_{j_{k}}$
with $$A^{'}_{j_{k}}:= \left\{ \bar{J}^{'}_{l}, \, \bar{j} -1 \geq l \geq j_{k}, \, and  \, |\bar{x}^{'}_{j_{k}} - \bar{x}^{'}_{l}| < |\bar{t}^{'}_{j_{k}}- \bar{t}^{'}_{l}| + 100R \right\}.$$
From Result \ref{Res:EstFirst} and the estimates above
\begin{align*}
A^{p+1} \card{\mathcal{P}}  & \gtrsim \sum_{j_{k} \in \mathcal{P}} \int_{J^{'}} \int_{\mathbb{R}^{3}}
\frac{|Iw(t,x)|^{p+1}}{|x-\bar{x}'_{j_{k}}|} \, dx \, dt \\
& \gtrsim  A^{-\alpha_{11}} \sum_{j_{k} \in \mathcal{P}} \int_{J^{'}} \int_{|x-\bar{x}^{'}_{j_{k}}| \leq |t-\bar{t}^{'}_{j_{k}}| + 1000R}
\frac{|I w(t,x)|^{p+1}} {1+ |t-\bar{t}_{j_{k}}^{'}| } \, dx \, dt \\
& \gtrsim A^{-\alpha_{11} - \alpha_{13}} \sum_{j_{k} \in \mathcal{P}} \sum_{l \in A^{'}_{j_{k}}} \int_{K^{'}_{l}} \frac{1}{1+ t} \, dt \\
& \geq c  A^{-\alpha_{11} - \alpha_{12} -\alpha_{13}} \sum_{j=1}^{\bar{j} -1} \frac{1}{1+ \bar{t}^{'}_{j}},
\end{align*}
where we used at the second line the elementary inequality $$\frac{1}{1000 R + |t-\bar{t}^{'}_{j_{k}}|} \gtrsim \frac{1}{R} \frac{1}{|t-\bar{t}^{'}_{j_{k}}| +1}.$$
Then it suffices to estimate $\card{\mathcal{P}}$. Let $j_{k_{max}}:= \max_{j_{k} \in \mathcal{P}} j_{k}$. By applying - $ \card{\mathcal{P}}$ -
times Result \ref{Res:IncreaseRate2}, by the construction of $\mathcal{P}$ and by (\ref{Eqn:LowerPotNrj})
\begin{align*}
A^{p+1} & \gtrsim E \left( I w(\bar{t}^{'}_{j_{k_{max}}}),
\bigcup_{j_{k} \in \mathcal{P}} B(\bar{x}^{'}_{j_{k}}, R + |\bar{t}^{'}_{j_{k_{max}}} - \bar{t}^{'}_{j_{k}}| ) \right)\\
& \geq \sum_{j_{k} \in \mathcal{P}} E \left( Iw(\bar{t}^{'}_{j_{k}}), B(\bar{x}^{'}_{j_{k}},R) \right) -
\frac{ c_{8} A^{-\alpha_{13}} \card{(\mathcal{P})} }{1000} \\
& \gtrsim  \card{(\mathcal{P})}  A^{- \alpha_{13}} \, .
\end{align*}
Hence, (\ref{Eqn:RepartJj}) follows.


\subsection{ Proof of Lemma \ref{lem:StrVsConcNrj}} $ $\\

Partitionning $J^{'}$ into the subintervals $(\bar{J}^{'}_{j})_{ 1 \leq j \leq \bar{j}}$ that were defined in Lemma \ref{lem:UseMorEst}, we see by  (\ref{Eqn:InducNrjw}) and
Proposition \ref{prop:LocalBd} that

\begin{equation}
Z(\bar{J}^{'}_{j},w)  \lesssim  A^{\frac{p+1}{2}}.
\label{Eqn:EstZJj}
\end{equation}
By (\ref{Eqn:ConcNrjMol}), there exists $(\bar{t}^{'}_{j},\bar{x}^{'}_{j}) \in \bar{J}^{'}_{j} \times \mathbb{R}^{3}$ such that
\begin{equation}
E \left( Iw(\bar{t}^{'}_{j}), B(\bar{x}^{'}_{j},R) \right)   \geq \frac{ c_{8} A^{-\alpha_{13}}}{p+1}
\end{equation}
Hence, from Result \ref{Res:IncreaseRate2} (see Appendix A), we see that there exist two constants $c_{9}$ and $\alpha_{16}$ such
that

\begin{equation}
E \left( I w(t), B(\bar{x}^{'}_{j}, R + |t -\bar{t}^{'}_{j}|)   \right)  \geq   c_{9} A^{-\alpha_{16}}, \, t \in \bar{J}^{'}_{j}.
\nonumber
\end{equation}
We further chop out each subinterval $\bar{J}^{'}_{j}$ into smaller subintervals
$(\bar{J}^{'}_{j,k}= [\bar{t}^{'}_{j,k}, \bar{t}^{'}_{j,k+1}])_{ k \in \mathbb{Z}}$ with
$\bar{t}^{'}_{j,0}:= \bar{t}_{j}^{'}$ and such that
$\tilde{Z}(\bar{J}^{'}_{j,k},w) \leq c_{9}  A^{- \alpha_{16}}$ while $\tilde{Z}(\bar{J}^{'}_{j,k},w) \sim  A^{- \alpha_{16}}$, except
maybe the last interval. Notice that by (\ref{Eqn:EstZJj}), there exists $\alpha$ such that
\begin{equation}
\forall j: \, \card  \left( ( \bar{J}^{'}_{j,k} )_{k \in \mathbb{Z}} \right)  \lesssim  A^{\alpha}.
\label{Eqn:UpBdK}
\end{equation}
From Result \ref{Res:EstFirst} we see that given $\epsilon > 0$, there exist $C >0$
and $\tilde{t}^{\hphantom{'}'}_{j,k}$ such that

\begin{equation}
\begin{array}{ll}
k \geq 0 : & \tilde{t}^{\hphantom{'}'}_{j,k} \in [\bar{t}^{'}_{j,k}, \bar{t}^{'}_{j,k} + C \langle \bar{t}^{'}_{j,k} - \bar{t}^{'}_{j} \rangle] \\
k < 0 : & \tilde{t}^{\hphantom{'}'}_{j,k} \in [\bar{t}^{'}_{j,k+1} - C \langle \bar{t}^{'}_{j,k+1} - \bar{t}^{'}_{j} \rangle,  \bar{t}^{'}_{j,k +1} ],
\end{array}
\nonumber
\end{equation}

\begin{equation}
\log{(C)}  \lesssim A^{\frac{(p+1)^{2}}{2}} \epsilon^{-p-1}, \qquad \left\| \frac{ I w (\tilde{t}^{\hphantom{'}'}_{j,k})}{ \langle x - \bar{x}^{'}_{j} \rangle} \right\|_{L^{2}}    \leq \epsilon.
\label{Eqn:BoundLogC0}
\end{equation}
\\ \\
Next, we claim that there exists $(k_{0},j_{0})$ and $M^{'}:= 1000  \langle C \rangle M$ such that
\begin{equation}
\begin{array}{ll}
k_{0} \geq 0: & | \bar{t}^{'}_{j_{0},k_{0}+1} - \bar{t}^{'}_{j_{0},k_{0}}|  \geq M^{'} \left( R + | \bar{t}^{'}_{j_{0},k_{0}}- \bar{t}^{'}_{j_{0}}| \right), \\ k_{0} < 0 : &  | \bar{t}^{'}_{j_{0},k_{0}+1} - \bar{t}^{'}_{j_{0},k_{0}}|  \geq M^{'} \left( R + | \bar{t}^{'}_{j_{0},k_{0}+1}- \bar{t}^{'}_{j_{0}}| \right).
\end{array}
\label{Eqn:DefMprime}
\end{equation}
If not, this implies, by simple induction on $k$, that (say)
$|\bar{t}^{'}_{j,k} - \bar{t}^{'}_{j}| \lesssim (2 M^{'})^{|k|+1} R$ for all $j$ and therefore, by (\ref{Eqn:UpBdK}), we see that there exist
two positive constants $\alpha$ and $C'$ such that $|\bar{J}^{'}_{j}| \leq (M^{'})^{ C' A^{\alpha}}$ for all $j$ \footnote{We allow the value of $\alpha$
 and $C'$ to increase in the sequel so that all the estimates
hold down to the end of the section.}. But, by Lemma \ref{lem:UseMorEst}, this imply that

\begin{equation*}
\frac{\log{\left( 1 + (\bar{j}-1) (M^{'})^{ C' A^{\alpha}} \right)}} { (M^{'})^{ C' A^{\alpha}} }  \leq \sum_{j=1}^{\bar{j}-1}
\frac{1}{1+ j  (M^{'})^{ C' A^{\alpha}}  }  \lesssim A^{\alpha_{14}} \, .
\end{equation*}
Therefore $ \log{(\bar{j})}  \leq (M^{'})^{ C' A^{\alpha}} $  and, combining this
inequality with (\ref{Eqn:BoundLogC0}), this yields a contradiction with (\ref{Eqn:UpperBdIu}). \\ \\
Assume that $k_{0} \geq 0$. Then
\begin{align*}
\bar{t}^{'}_{j_{0},k_{0}+1} -\tilde{t}^{\hphantom{'}'}_{j_{0},k_{0}} & \geq \bar{t}^{'}_{j_{0},k_{0}+1} - \bar{t}^{'}_{j_{0},k_{0}} - C
\langle \bar{t}^{'}_{j_{0},k_{0}}
- \bar{t}^{'}_{j_{0}} \rangle \\
& \geq M^{'} \left(  R + | \bar{t}^{'}_{j_{0},k_{0}} - \bar{t}^{'}_{j_{0}}  | \right)   - C \langle \bar{t}^{'}_{j_{0},k_{0}} - \bar{t}^{'}_{j_{0}} \rangle \\
& \geq M \left( R + |\tilde{t}^{\hphantom{'}'}_{j_{0},k_{0}} - \bar{t}^{'}_{j_{0}} | \right)
\end{align*}
Hence, choosing $\bar{x}:= \bar{x}^{'}_{j_{0}}$, $R^{'}:= R + | \tilde{t}^{\hphantom{'}'}_{j_{0},k_{0}} -\bar{t}^{'}_{j_{0}} | $, $S:= \tilde{t}^{\hphantom{'}'}_{j_{0},k_{0}}$ and $T:= \bar{t}^{'}_{j_{0},k_{0}+1}$, we have (\ref{Eqn:BoundZpr}), (\ref{Eqn:LowerBoundST}) and (\ref{Eqn:UpBdWeight}). The reader is invited to check
that if $k_{0} < 0$, then a similar argument shows that the same estimates hold if  $R^{'}:= R + | \tilde{t}^{\hphantom{'}'}_{j_{0},k_{0}+1} -\bar{t}^{'}_{j_{0}} | $, $S:= \tilde{t}^{\hphantom{'}'}_{j_{0},k_{0}+1}$ and $T:= \bar{t}^{'}_{j_{0},k_{0}}$.


\section{Proof of perturbation argument}
\label{Sec:ProofPerturb}
\noindent In this section we prove Proposition \ref{prop:PerturbArg}. We may assume without loss of generality that
$J^{''}=(T,b)$. The proof is made of several steps: \\
\\
Step 1 . {\it Bound of $Z([T,b],w)$ and  $Z([T,b],w')$  }. \\
\\
We divide $[T,b]$ into subintervals $K$ such that $ \| I w \|_{L_{t}^{2(p-1)-} L_{x}^{2(p-1)+} (K)} = c_{2} $, except maybe the last
one. Proposition \ref{prop:LocalBd} yields $Z(K,w) \lesssim A^{\frac{p+1}{2}}$. Iterating over $K$ and using (\ref{Eqn:InducStr1})

\begin{equation}
Z([T,b],w)  \lesssim  \langle L_{w} \rangle^{2(p-1)-} A^{\frac{p+1}{2}}.
\label{Eqn:BoundwlPerturb}
\end{equation}
A bound of $Z(w',[T,b])$ is obtained similarly as above:
\begin{equation}
Z([T,b],w')  \lesssim  \langle L_{w'} \rangle^{2(p-1)-} A^{\frac{p+1}{2}} \, .
\label{Eqn:BoundwlplusonePerturb}
\end{equation}

Step 2.  {\it Decomposition}. \\
\\
Let $\Gamma:= w - w^{'} - v$ and $K^{'}:=[t',t''] \subset [T,b]$. A simple computation shows that

\begin{equation*}
\partial_{tt} I \Gamma - \triangle I \Gamma + I \Gamma  = (I F(w^{'}) - F(Iw^{'})) + ( F(Iw) - I F(w)) + ( F(Iw^{'}) - F(Iw)).
\nonumber
\end{equation*}
We decompose

\begin{align*}
 I \Gamma(t)  = I \Gamma^{K'}_{l}(t) + X^{K'}_{F(Iw') - I F(w')}(t) + X^{K'}_{I F(w) - F(Iw)}(t) + X^{K'}_{F(Iw)-F(Iw')}(t),
\end{align*}
where the $X$ numbers are defined in Section \ref{Sec:Notation} and

\begin{align*}
\begin{array}{ll}
\Gamma_{l}^{K'}(t) & := \cos{ ((t-t^{'})\langle D \rangle)} \Gamma(t^{'}) + \frac{\sin{ \left((t -t^{'}) \langle D  \rangle \right)}}{\langle D \rangle} \partial_{t} \Gamma(t^{'}).
\end{array}
\nonumber
\end{align*}
We use the following notation: Given a function f, let
\begin{equation}
\bar{Z}(K',f)  := \sup_{(q,r)-\frac{1}{2} \,  wave \, adm} \left\| \langle D \rangle^{s_{c} - \frac{1}{2}} I f \right\|_{L_{t}^{q} L_{x}^{r} (K')}.
\label{Eqn:DfnZbqrperturb}
\end{equation}
\\
Step 3. {\it Short-time perturbation argument.}
\begin{lem}
There exist $ 0 < \theta:=\theta(p) < 1$ and $c_{10}$ such that if $c \leq c_{10}$, \\

\begin{align}
\| I w' \|_{L_{t}^{2(p-1)-} L_{x}^{2(p-1)+}(K')}, \, \| \langle D  \rangle ^{s_{c} -\frac{1}{2}} I w' \|_{L_{t}^{4+} L_{x}^{4-}(K')}  \leq c,
\label{Eqn:Hyp0} \\
\bar{Z}(K',\Gamma_{l}^{K'}) & \leq c, \, and
\label{Eqn:Hyp1} \\
\| I v \|_{L_{t}^{2(p-1)-} L_{x}^{2(p-1)+} (K') } & \leq  c,
\label{Eqn:Hyp2} \\
\nonumber
\end{align}
then
\begin{equation}
\begin{array}{l}
\bar{Z}(K', I^{-1} X^{K'}_{F(Iw')-IF(w')}) +  \bar{Z}(K', I^{-1} X^{K'}_{IF(w)-F(Iw)} ) +
\bar{Z}(K', I^{-1} X^{K'}_{F(Iw)-F(Iw')} )  \\ \\
\lesssim \max  \left( \frac{\max^{2p(p-1)-}(\langle L_{w} \rangle, \langle L_{w'} \rangle) A^{\frac{p(p+1)}{2}}}{N^{(1-s_{c})-}}, \,  c \right).
\end{array}
\label{Eqn:Concl}
\end{equation}
\label{lem:ShortTimePerturb}
\end{lem}

\begin{proof}

We have

\begin{equation}
\begin{array}{l}
\bar{Z}(K',\Gamma ) \leq \bar{Z}(K',\Gamma^{K'}_{l}) + \bar{Z} (K',I^{-1} X^{K'}_{F(Iw')-IF(w')})
+ \bar{Z}(K',I^{-1} X^{K'}_{IF(w)-F(Iw)}) + \bar{Z}(K',I^{-1} X^{K'}_{F(Iw) - F(Iw')}) \cdot
\end{array}
\nonumber
\end{equation}
We first estimate $ \bar{Z}(K', I^{-1} X_{F(Iw)-F(Iw')}^{K'} )$. \\
By interpolation (see points $A$, $B$, and $C$ on Figure \ref{Fig:StrichGraph}) there exist $\theta:= \theta(p) > 0$, $m < \frac{1}{2}$, and $(\bar{q},\bar{r})$ $m$-wave admissible such that

\begin{equation}
\begin{array}{l}
\| \langle D  \rangle^{\frac{1}{2}} I v \|_{L_{t}^{4+} L_{x}^{4-} (K')} \lesssim
\| \langle D \rangle^{1-s_{c}} I v \|^{\theta}_{L_{t}^{2(p-1)-} L_{x}^{2(p-1)+} (K')}
 \| \langle D \rangle^{1-m}  I v  \|^{1- \theta}_{L_{t}^{\bar{q}} L_{x}^{\bar{r}} (K')}
\end{array}
\label{Eqn:Interp}
\end{equation}
Notice that, in view of (\ref{Eqn:ConservkinetIv}) and (\ref{Eqn:StrNlkg}), we have
$\| \langle D \rangle^{1-m} I v \|_{L_{t}^{\bar{q}} L_{x}^{\bar{r}} (K')} \lesssim  \| \langle D \rangle I v  (t') \|_{L^{2}} \lesssim 1 $. Hence
\begin{align*}
\| \langle D  \rangle^{s_{c} - \frac{1}{2}} I v \|_{L_{t}^{4+} L_{x}^{4-} (K^{'})} \lesssim
\| I v \|^{\theta}_{L_{t}^{2(p-1)-} L_{x}^{2(p-1)+} (K') } \lesssim  c^{\theta}.
\end{align*}
Therefore, using again (\ref{Eqn:StrNlkg})
\begin{equation}
\begin{array}{l}
\| \langle D \rangle^{ s_{c} - \frac{1}{2} }  X_{F(Iw)-F(Iw')}^{K'} \|_{L_{t}^{q} L_{x}^{r} (K')}  \\
\\
\lesssim \| \langle D \rangle^{ s_{c} - \frac{1}{2} }
(F(Iw) - F(Iw') \|_{L_{t}^{\frac{4}{3}} L_{x}^{\frac{4}{3}} (K')} \\
\\
\lesssim
\left(
\| I w \|^{p-1}_{ L_{t}^{2(p-1)-} L_{x}^{2(p-1)+} (K') }  + \| I w^{'} \|^{p-1}_{L_{t}^{2(p-1)-} L_{x}^{2(p-1)+} (K') }
\right) \left\| \langle D \rangle^{s_{c} -\frac{1}{2}} I (w- w') \right\|_{L_{t}^{4+} L_{x}^{4-} (K')} \\
+ \left(  \| I w \|^{p-2}_{L_{t}^{2(p-1)-} L_{x}^{2(p-1)+} (K')} +  \| I w' \|^{p-2}_{L_{t}^{2(p-1)-} L_{x}^{2(p-1)+} (K')} \right)
\left(  \| \langle D \rangle^{s_{c} -\frac{1}{2}} I w \|_{L_{t}^{4+} L_{x}^{4-}(K')} + \| \langle D \rangle^{s_{c} -\frac{1}{2}} I w' \|_{L_{t}^{4+} L_{x}^{4-}(K')} \right) \\
\| I(w-w') \|_{L_{t}^{2(p-1)-} L_{x}^{2(p-1)+} (K')} \\
\\
\lesssim
 \left(
\| I w \|^{p-1}_{ L_{t}^{2(p-1)-} L_{x}^{2(p-1)+} (K') }  + \| I w^{'} \|^{p-1}_{L_{t}^{2(p-1)-} L_{x}^{2(p-1)+} (K') }
\right)  \times \left( \| \langle D \rangle^{ s_{c} - \frac{1}{2} } I v \|_{ L_{t}^{4+} L_{x}^{4-} (K')}
 + \| \langle D \rangle^{ s_{c} - \frac{1}{2}} I \Gamma \|_{ L_{t}^{4+} L_{x}^{4- } (K') } \right) \\
+ \left(  \| I w \|^{p-2}_{L_{t}^{2(p-1)-} L_{x}^{2(p-1)+} (K')} +  \| I w' \|^{p-2}_{L_{t}^{2(p-1)-} L_{x}^{2(p-1)+} (K')} \right)
\left(  \| \langle D \rangle^{s_{c} -\frac{1}{2}} I w \|_{L_{t}^{4+} L_{x}^{4-}(K')} + \| \langle D \rangle^{s_{c} -\frac{1}{2}} I w' \|_{L_{t}^{4+} L_{x}^{4-}(K')} \right) \\
\left( \| I v \|_{L_{t}^{2(p-1)-} L_{x}^{2(p-1)+} (K')} + \| I \Gamma \|_{L_{t}^{2(p-1)-} L_{x}^{2(p-1)+} (K')} \right) \\
\\
\lesssim   \left( \bar{Z}^{p-1}(K',\Gamma)  + c^{p-1} \right)
\left( c^{\theta} + \bar{Z}(K',\Gamma) \right)
 + \left( \bar{Z}^{p-2}(K',\Gamma) + c^{p-2} \right) \left( \bar{Z}(K',\Gamma) + c^{\theta} \right)
\left( \bar{Z}(K',\Gamma) + c \right),
\end{array}
\nonumber
\end{equation}
substituting $w$  for $w' + \Gamma + v$ and using the Sobolev embedding at the last line, i.e
\begin{align}
\| I \Gamma \|_{ L_{t}^{2(p-1)-} L_{x}^{2(p-1)+} (K')} \lesssim  \|  \langle D \rangle^{s_{c}- \frac{1}{2}} I \Gamma \|_{
L_{t}^{2(p-1)-} L_{x}^{\frac{6(p-1)}{2p-3} +} (K') }  \lesssim  \bar{Z}(K',\Gamma).
\label{Eqn:SobEmbedSpec}
\end{align}
Hence, collecting all these estimates, and using Result \ref{Eqn:EstZb} (see Appendix A), we get (\ref{Eqn:Concl}) from a continuity argument.

\end{proof}

\noindent Step 4. {\it Long-time perturbation argument.} \\
\\
We divide $[T,b]$ into subintervals $(K'_{q}=[t'_{q}, t'_{q+1}])_{1 \leq q \leq Q}$ such that
$\| I w' \|_{L_{t}^{2(p-1)-} L_{x}^{2(p-1)+} (K'_{q})} = c$ or $\| \langle D \rangle^{s_{c} -\frac{1}{2}} I w' \|_{L_{t}^{4+} L_{x}^{4-} (K'_{q})} = c$,
 while $\max ( \| I w' \|_{L_{t}^{2(p-1)-} L_{x}^{2(p-1)+} (K'_{q})}, \| \langle D \rangle^{s_{c} -\frac{1}{2}} I w' \|_{L_{t}^{4+} L_{x}^{4-} (K'_{q})} ) \leq c$, except maybe the last one.
We can use the short-time perturbation argument on $K^{'}:=J'_{q}$ as long as
(\ref{Eqn:Hyp1}) and (\ref{Eqn:Hyp2}) hold (with $c \leq c_{10}$ ). But since
$\Gamma_{l}(T)=0$ and

\begin{equation*}
\bar{Z} (K'_{q'}, \Gamma_{l}^{K'_{q'}} )
\lesssim
\bar{Z} ( K'_{q'}, \Gamma_{l}^{[T,b]} )
+ \sum_{q=1}^{q'-1} \bar{Z}(K'_{q}, I^{-1} X^{K'_{q}}_{F(Iw')-IF(w')})
+ \bar{Z}(K'_{q}, I^{-1} X^{K'_{q}}_{IF(w)-F(Iw)})
+  \bar{Z}(K'_{q}, I^{-1} X^{K'_{q}}_{F(Iw)-F(Iw')} )
\nonumber
\end{equation*}
we easily see, by iteration, that there exists a positive constant $C$ such that

\begin{equation*}
\bar{Z}(K'_{q'}, \Gamma_{l}^{K'_{q'}} ) \leq C^{q'}
\max \left( k, \frac{ A^{\frac{p(p+1)}{2}} \max^{2p (p-1)-}( \langle L_{w} \rangle , \langle L_{w'} \rangle)}{N^{(1-s_{c})-}} \right).
\end{equation*}
Hence we see from (\ref{Eqn:BoundwlplusonePerturb}), (\ref{Eqn:Cond4}) and (\ref{Eqn:Boundk}), we see that (\ref{Eqn:Hyp1}) and (\ref{Eqn:Hyp2}) hold, by choosing $C_{3}$, $C_{4}$, $\alpha_{6}$, $ \alpha_{7} $ (resp. $c_{4}$) large enough (resp. small enough). \\
By summation over $q^{'}$ and (\ref{Eqn:SobEmbedSpec}) (with $K'$ substituted for $[T,b]$), we see that there exists a positive constant $\alpha$
such that $\| I \Gamma \|_{L_{t}^{2(p-1)-} L_{x}^{2(p-1)+} ([T,b])} \lesssim  \left( \langle L_{w'} \rangle A \right)^{\alpha}$. Hence (\ref{Eqn:Cond5Res}) holds.

\section{Proof of small mollified energy theory}
\label{Sec:EstSmallNrj}

\noindent In this section we prove Proposition \ref{prop:EstSmallNrj}.
The proof is made of two steps: \\
\\
\noindent \textit{Control of $\| I w \|_{L_{t}^{2(p-1)-} L_{x}^{2(p-1)+} (\mathbb{R})}$}. \\
\\
By (\ref{Eqn:SmallMolNrjTime}) and (\ref{Eqn:Comput}) we realize that
\begin{equation}
\begin{array}{ll}
Z_{s_{c},s}(\mathbb{R},w) & \lesssim  E^{\frac{1}{2}}(Iw(\tilde{t}))  +   \| \langle D \rangle^{1-s_{c}} I w  \|_{L_{t}^{\frac{2}{s_{c}} +} L_{x}^{\frac{2}{1-s_{c}}-} (\mathbb{R})}
\times
\left(
\begin{array}{l}
 \| Iw \|^{p-1}_{L_{t}^{2(p-1)-} L_{x}^{2(p-1)+} (\mathbb{R})}  \\
+ \frac{ \| \langle D \rangle^{1-s_{c}} I w  \|^{p-1}_{L_{t}^{2(p-1)-} L_{x}^{2(p-1)+} (\mathbb{R}) } } {N^{\frac{5-p}{2}-}}
\end{array}
\right)
 \\
& \lesssim E^{\frac{1}{2}}(Iw(\tilde{t})) + \left( Z^{p}_{s_{c},s}(\mathbb{R},w) + \frac{Z^{p}_{s_{c},s} (\mathbb{R},w)}{N^{\frac{5-p}{2}-}}  \right),
\end{array}
\nonumber
\end{equation}
where we used the Sobolev embedding, that is
$$ \| I w \|_{L_{t}^{2(p-1)-} L_{x}^{2(p-1)+} (\mathbb{R})}  \lesssim \| \langle D \rangle^{1- s_{c}} I w \|_{L_{t}^{2(p-1)-} L_{x}^{\frac{6(p-1)}{2p-3}+ }(\mathbb{R})}. $$
Therefore we see by a continuity argument that  $Z_{s_{c},s}(\mathbb{R},w) \lesssim E^{\frac{1}{2}}(Iw(\tilde{t}))$ and, consequently, (\ref{Eqn:ControlIw}) holds.\\
\\
\noindent  \textit{Control of $\sup_{t \in \mathbb{R}} E(Iw(t))$}. \\
\\

Let $T > 0$. From (\ref{Eqn:ControlIw}) one may divide $[\tilde{t}, \tilde{t} + T]$ and $[\tilde{t} - T, \tilde{t}]$ into subintervals $J$ such that
$\| I w \|_{L_{t}^{2(p-1)-} L_{x}^{2(p-1)+} (J)} =c_{2}$, except maybe the last one. By choosing $c_{5}$ (resp. $\alpha_{9}$) small enough
(resp. large enough) in (\ref{Eqn:InducStr3}), we see that we may apply Proposition \ref{prop:LocalBd} and Proposition \ref{Prop:Acl} on each $J$ and we get after iteration that

\begin{align*}
|E(Iw(\tilde{t} + T)) - E(Iw(\tilde{t}))|, \, |E(Iw(\tilde{t} - T)) - E(Iw(\tilde{t}))|   \lesssim \frac{E^{\frac{p+1}{2}} (I w(\tilde{t}))}{N^{\frac{5-p}{2}-}},
\end{align*}
assuming that $\sup_{t \in [\tilde{t} - T, \tilde{t} + T] } E(Iw(t)) \lesssim E(Iw(\tilde{t})) $. But, since again $N$ satisfies
(\ref{Eqn:InducStr3}) with $c_{5}$ (resp. $\alpha_{9}$) small enough (resp. large enough), we see that not only this estimate holds but also (\ref{Eqn:ControlSupNrjSmall}).

\section{Appendix A: Estimates involving commutators}

\subsection{ Estimates involving commutators in Section \ref{Sec:SepLocMol}}

We prove all the estimates involving commutators that appear in Section \ref{Sec:SepLocMol}.

\subsubsection{Result \ref{Eqn:BoundX11}}

\begin{res}
Let $X_{1,1}$, $X_{1,2}$ be defined in (\ref{Eqn:DefX1112}). Then

\begin{equation}
|X_{1,1}|, |X_{1,2}|  \leq \frac{c_{9} A^{- \alpha_{16}}}{1000}
\nonumber
\end{equation}
\label{Eqn:BoundX11}
\end{res}

\begin{proof}

In view of (\ref{Eqn:InducNrjw}) and (\ref{Eqn:IneqNrjS})

\begin{equation}
\| \partial_{t} I \bar{w} \|_{L_{t}^{\infty} L_{x}^{2}([S,T])}  \lesssim  A^{\frac{p+1}{2}}.
\label{Eqn:EstDerivBarw}
\end{equation}
We first estimate $X_{1,2}$. From (\ref{Eqn:InducStr2}), (\ref{Eqn:BoundZpr}), (\ref{Eqn:EstDerivBarw}), and \\
$|F(w) - F(Iw)| \lesssim \max{ ( |Iw|^{p-1}, |w|^{p-1})} | I w -w|$,
\begin{align*}
|X_{1,2}|  & \lesssim \| \partial_{t} I \bar{w} \|_{L_{t}^{\infty} L_{x}^{2}([S,T])} \times \left(
\begin{array}{l}
\| P_{<<N} w \|^{p-1}_{L_{t}^{\frac{4(p-1)}{7-p}} L_{x}^{\frac{4(p-1)}{p-3} } ([S,T]) } \| P_{\gtrsim N} w
\|_{L_{t}^{\frac{4}{p-3}}
L_{x}^{\frac{4}{5-p}} ([S,T]) } \\
+ \| P_{\gtrsim N} w \|^{p-1}_{L_{t}^{p} L_{x}^{2p} ([S,T])} \| P_{\gtrsim N} w \|_{L_{t}^{p} L_{x}^{2p} ([S,T])}
\end{array}
\right) \\
& \lesssim \frac{1}{N^{\frac{5-p}{2}-}} \| \partial_{t} I \bar{w} \|_{L_{t}^{\infty} L_{x}^{2}([S,T])} \times \left(
\begin{array}{l}
\| \langle D \rangle^{1- 1} I w \|^{p-1}_{L_{t}^{\frac{4(p-1)}{7-p} } L_{x}^{\frac{4(p-1)}{p-3} } ([S,T]) } \\ \times
\| \langle D \rangle^{1-  \frac{p-3}{2}}  I w  \|_{L_{t}^{\frac{4}{p-3}} L_{x}^{\frac{4}{5-p}} ([S,T]) } \\
+ \| \langle  D \rangle^{1- \frac{3p-5}{2p}} I w  \|^{p}_{L_{t}^{p} L_{x}^{2p} ([S,T])}
\end{array}
\right) \\
& \leq  \frac{c_{9} A^{-\alpha_{16}}}{1000},
\end{align*}
where at the last line we choose $\alpha_{1}$ (resp. $c_{1}$) large enough (resp. small enough) in (\ref{Eqn:InducStr2}). \\
We turn to $X_{1,1}$. We use an argument in \cite{triroyradnlkg}: For low frequencies we use the smoothness of $F$ ($F$ is $C^{1}$) and for high
frequencies, we use the regularity of $w$ (in $H^{s}$). Indeed, we have
\begin{align}
F (w) &:= F \left( P_{\ll N} w + P_{\gtrsim N} w  \right)  \nonumber \\
&= F \left( P_{\ll N} w \right) + \left( \int_{0}^{1} |P_{\ll N} w + y P_{\gtrsim N} w|^{p-1} \, dy \right) P_{\gtrsim N} w \nonumber \\
&+ \left( \int_{0}^{1} \frac{ P_{\ll N} w + y P_{\gtrsim N} w }{  _{ \quad  \overline{P_{\ll N} w + y P_{\gtrsim N} w} \quad  } } |P_{\ll N} w + y
P_{\gtrsim N} w|^{p-1} \, dy \right) \, \overline{P_{\gtrsim N} w} \, .
\label{Eqn:DecompFLowHigh}
\end{align}
Therefore, we estimate
\begin{align*}
|X_{1,1}| & \lesssim \| \partial_{t} I \bar{w} \|_{L_{t}^{\infty} L_{x}^{2}([S,T])} (X_{1,1,1} + X_{1,1,2} + X_{1,1,3} ) \\
& \lesssim A^{\frac{p+1}{2}}  (X_{1,1,1} + X_{1,1,2} + X_{1,1,3} )
\end{align*}
with
\begin{align*}
 X_{1,1,1} & :=  \| P_{\gtrsim N} F(P_{\ll N} w) \|_{L_{t}^{1} L_{x}^{2}([S,T])}   \\
X_{1,1,2} & := \| | P_{\ll N} w |^{p-1} P_{\gtrsim N} w \|_{L_{t}^{1} L_{x}^{2}([S,T])} \\
X_{1,1,3} & :=  \| P_{\gtrsim N} w \|^{p}_{L_{t}^{p} L_{x}^{2p} ([S,T])}.
\end{align*}
We further estimate
\begin{align*}
X_{1,1,1} & \lesssim \frac{1}{N^{1-}} \| \nabla F(P_{\ll N} w) \|_{L_{t}^{1} L_{x}^{2}([S,T])} \\
& \lesssim \frac{1}{N ^{1-}} \| P_{\ll N} w \|^{p-1}_{L_{t}^{\frac{4(p-1)}{7-p}} L_{x}^{\frac{4(p-1)}{p-3}} ([S,T])} \| \nabla P_{\ll N} w
\|_{L_{t}^{\frac{4}{p-3}} L_{x}^{\frac{4}{5-p}} ([S,T])} \\
& \lesssim \frac{1}{N^{\frac{5-p}{2}-}} \| \langle D \rangle^{1-1} I w  \|^{p-1}_{L_{t}^{\frac{4(p-1)}{7-p} } L_{x}^{\frac{4(p-1)}{p-3}}
([S,T]) } \| \langle D \rangle^{1- \frac{p-3}{2}} I w \|_{L_{t}^{\frac{4}{p-3}} L_{x}^{\frac{4}{5-p} }([S,T])}  \\
& \lesssim \frac{ A^{- p \alpha_{16}}}{N^{\frac{5-p}{2}-}} \, ,
\end{align*}

\begin{align*}
X_{1,1,2} & \lesssim \left\| |P_{\ll N} w|^{p-1} P_{\gtrsim N} w \right\|_{L_{t}^{1} L_{x}^{2} ([S,T])} \\
& \lesssim \frac{1}{N^{\frac{5-p}{2}-}} \| \langle D \rangle^{1-1} I w \|^{p-1}_{L_{t}^{\frac{4(p-1)}{7-p}} L_{x}^{\frac{4(p-1)}{p-3}}
([S,T])}
\| \langle D \rangle^{1- \frac{p-3}{2} } I w  \|_{L_{t}^{\frac{4}{p-3}} L_{x}^{\frac{4}{5-p}} ([S,T]) } \\
& \lesssim \frac{ A^{- p \alpha_{16}}}{N^{\frac{5-p}{2}-}} \, ,
\end{align*}
and
\begin{align*}
X_{1,1,3} & \lesssim \| P_{\gtrsim N} w \|^{p}_{L_{t}^{p} L_{x}^{2p} ([S,T])} \\
& \lesssim \frac{ \| \langle D \rangle^{1- \frac{3p-5}{2p}} I w \|^{p}_{L_{t}^{p} L_{x}^{2p}([S,T])}} {N^{\frac{5-p}{2}-}} \\
& \lesssim \frac{ A^{ - p \alpha_{16}}}{N^{\frac{5-p}{2}-}} \, .
\end{align*}
Combining these estimates and using again (\ref{Eqn:InducStr2}) with $\alpha_{1}$ (resp. $c_{1}$) large enough (resp. small enough), we obtain
$|X_{1,1}| \leq  \frac{c_{9} A^{-\alpha_{16}}}{1000} $.

\end{proof}

\subsubsection{Result \ref{Res:EstFirst}}

\begin{res}
Let $(\tilde{t},\tilde{x}) \in \mathbb{R} \times \mathbb{R}^{3}$. Then

\begin{equation}
\int_{J} \int_{\mathbb{R}^{3}} \frac{|Iw(t,x)|^{p+1}}{|x-\tilde{x}|} \, dx \, dt  \lesssim A^{p+1}  + \frac{ \langle L_{w} \rangle^{2(p-1)-} A^{\frac{(p+1)^{2}}{2}} }
{N^{\frac{5-p}{2}-}}   \lesssim  A^{p+1},
\label{Eqn:EstLongEst}
\end{equation}
and

\begin{equation}
\int_{J} \left( \int_{\mathbb{R}^{3}} \frac{|Iw|^{2}}{ \langle x - \tilde{x} \rangle^{2}} \, dx \right)^{\frac{p+1}{2}} \frac{1}{ \langle t -\tilde{t} \rangle} \, dt \lesssim A^{\frac{(p+1)^{2}}{2}} \cdot
\label{Eqn:LongTimeHardy}
\end{equation}

\label{Res:EstFirst}
\end{res}

\begin{proof}

From (\ref{Eqn:InducStr1}), one can  chop $J$ into subintervals $K$ such that
$\| I w \|_{L_{t}^{2(p-1) -} L_{x}^{2(p-1) +}(K)} = c_{2}$, except maybe the last one.
From (\ref{Eqn:InducNrjw}), Proposition \ref{prop:LocalBd} and Proposition \ref{Prop:AlmMor}, we see, by iteration, that
\begin{equation}
\int_{J} \int_{\mathbb{R}^{3}} \frac{|Iw(t,x)|^{p+1}}{|x- \tilde{x}|} \, dx \, dt  \lesssim
 A^{p+1}  + \frac{ \langle L_{w} \rangle^{2(p-1)-} A^{\frac{(p+1)^{2}}{2}} }
{N^{\frac{5-p}{2}-}} .
\nonumber
\end{equation}
Choosing $\alpha_{1}$ (resp. $c_{1}$) large enough (resp. small enough) in (\ref{Eqn:InducStr2}), we get (\ref{Eqn:EstLongEst}). \\
Next we prove (\ref{Eqn:LongTimeHardy}). Following Lemma $5.3$, \cite{nakanash1}, we have
\begin{equation*}
\int_{J} \left( \int_{\mathbb{R}^{3}} \frac{|Iw|^{2}}{\langle x - \tilde{x} \rangle^{2}} \, dx \right)^{\frac{p+1}{2}} \frac{1}{\langle t - \tilde{t}
\rangle} \, dt  \leq X_{1} + X_{2}
\end{equation*}
where
\begin{align*}
X_{1} &:= \int_{J} \left( \int_{|t - \tilde{t}| > |x- \tilde{x}|}  \frac{|I w|^{2}}{\langle x -  \tilde{x} \rangle^{2}} \, dx  \right)^{\frac{p+1}{2}} \frac{1}{\langle t - \tilde{t} \rangle} \, dt,  \\
X_{2} &:= \int_{J} \left( \int_{|t- \tilde{t}| \leq |x- \tilde{x}|} \frac{|I w|^{2}}{\langle x - \tilde{x} \rangle^{2}} \, dx \right)^{\frac{p+1}{2}} \frac{1}{\langle
t - \tilde{t} \rangle} \, dt.
\end{align*}
By H\"older inequality and (\ref{Eqn:EstLongEst})
\begin{align*}
X_{1}   \leq \int_{J}  \int_{|t - \tilde{t}| > |x - \tilde{x}|}  \frac{|I w|^{p+1}}{\langle t - \tilde{t} \rangle} \, dx \, dt \lesssim A^{p+1}
\end{align*}
We also have
\begin{align*}
X_{2}\leq \int_{J}   \sup_{|x - \tilde{x}| \geq |t - \tilde{t}|} \left(  \frac{1}{ \langle  x - \tilde{x} \rangle^{2}} \right)^{\frac{p+1}{2}}
\left( \int_{|x- \tilde{x}| \geq |t - \tilde{t}|} |Iw|^{2}  \right)^{\frac{p+1}{2}} \frac{1}{\langle t - \tilde{t} \rangle} \, dt \lesssim A^{\frac{(p+1)^{2}}{2}}.
\end{align*}

\end{proof}

\subsubsection{Result \ref{Res:IncreaseRate2}}

\begin{res}
Let $x^{'}_{a} \in \mathbb{R}^{3}$ and let $t^{'}_{b} \geq t^{'}_{a}$. Then we have
\begin{equation}
E \left( Iw(t^{'}_{b}), B(x^{'}_{a}, R + t^{'}_{b} -t^{'}_{a}) \right)  \geq E \left( Iw(t^{'}_{a}), B(x^{'}_{a},R)  \right) - \frac{c_{8} A^{-\alpha_{13}}}{1000}  \, .
\label{Eqn:IncreaseRate2}
\end{equation}
\label{Res:IncreaseRate2}
\end{res}

\begin{proof}

Integrating the identity
\begin{align*}
\partial_{t} \big( & \frac{1}{2} |\partial_{t} I w|^{2} + \frac{1}{2} |\nabla I w|^{2} + \frac{|Iw|^{p+1}}{p+1} + \frac{|Iw|^{2}}{2} \big)
- \partial_{x_{i}} \big( \Re ( \overline{\partial_{t} I w} \partial_{x_{i}} I w )   \big) \\
&+ \Re \big( \overline{\partial_{t} I w} (I F(w) - F(Iw)) \big)=0,
\end{align*}
inside the truncated cone $M:= \{ (t,x), \, t \in (t^{'}_{a}, t^{'}_{b}), \,  t-t^{'}_{a} - R \geq |x- x^{'}_{a}|  \}$, we obtain
\begin{align} \label{Eqn:IncreaseRate}
E &\left( I w(t^{'}_{b}), B( x^{'}_{a}, R + t^{'}_{b} - t^{'}_{a}) \right) - E \left( I w(t^{'}_{a}), B(x^{'}_{a}, R) \right) \\ \nonumber
&= \frac{1}{ \sqrt{2}} \int_{\partial M} \frac{|Iw|^{2}}{2} + \frac{|Iw|^{p+1}}{p+1} \, d \sigma \\ \nonumber
&+ \frac{1}{\sqrt{2}} \int_{\partial M} \left| \frac{x-x^{'}_{a}}{|x-x^{'}_{a}|} \partial_{t} I w + \nabla I w   \right|^{2} \, d \sigma -
\int_{M} \Re \left( \overline{\partial_{t} I w} (I F(w) - F(Iw))   \right) \, dx dt \, .
\end{align}
The boundary terms are nonnegative. In order to deal with the last integral, we chop $J$ into subintervals $K$ such that
$\| I w \|_{L_{t}^{2(p-1) -} L_{x}^{2(p-1) +}(K)} = c_{2}$, except maybe the last one; then, from (\ref{Eqn:InducStr1}), Proposition \ref{Prop:Acl}, Proposition \ref{prop:LocalBd}, and iteration, we get

\begin{align*}
\left| \int_{M} \Re \left( \overline{\partial_{t} I w} (I F(w) - F(Iw))   \right) \, dx dt \right| & \lesssim \frac{ \langle L_{w} \rangle^{2(p-1)-} A^{\frac{(p+1)^{2}}{2}}}{ N^{\frac{5-p}{2}-}} \\
& \leq \frac{c_{8} A^{-\alpha_{13}}}{1000},
\end{align*}
choosing $\alpha_{1}$ (resp. $c_{1}$) large enough (resp. small enough) in (\ref{Eqn:InducStr2}).

\end{proof}






\subsection{Estimates involving commutators in Section \ref{Sec:ProofPerturb}}
We prove all the estimates involving commutators that appear in Section \ref{Sec:ProofPerturb}.

\subsubsection{Result \ref{Eqn:EstZb}}

\begin{res}
\begin{equation}
\bar{Z}(K', I^{-1} X_{IF(w) - F(Iw)}^{K'}), \bar{Z}(K', I^{-1} X_{F(Iw')-I F(w')}^{K'}) \lesssim \frac{\max^{2p(p-1)-}{( \langle L_{w} \rangle , \langle L_{w'} \rangle)} A^{\frac{p(p+1)}{2}}}{N^{(1-s_{c})-}}
\nonumber
\end{equation}
\label{Eqn:EstZb}
\end{res}

\begin{proof}

Step 1 . {\it Bound of  $\bar{Z}(K', I^{-1} X^{K'}_{IF(w)-F(Iw)} )$ and   $\bar{Z}(K', I^{-1} X^{J'}_{F(Iw')-IF(w')} )$}. \\
\\
We first estimate $\bar{Z}(K', I^{-1} X^{K'}_{IF(w)-F(Iw)} )$. \\
We write $X^{K'}_{IF(w)-F(Iw)} = X^{K'}_{IF(w)-F(w)} + X^{K'}_{F(w)-F(Iw)}$.
By (\ref{Eqn:BoundwlPerturb}) and (\ref{Eqn:StrNlkg}) (with $(q,r)$ defined in (\ref{Eqn:DfnZbqrperturb}))
\begin{equation}
\begin{array}{ll}
\| \langle D  \rangle^{s_{c} - \frac{1}{2}} X^{K'}_{F(w)-F(Iw)} \|_{ L_{t}^{q} L_{x}^{r} (K^{'}) } &  \lesssim
\| \langle D \rangle ^{s_{c} - \frac{1}{2}} ( F(w)-F(Iw)) \| _{L_{t}^{\frac{4}{3}}  L_{x}^{\frac{4}{3}}  (K^{'}) } \\
\\
& \lesssim \| w \|^{p-1}_{L_{t}^{2(p-1)} L_{x}^{2(p-1)} ([T,b]) } \| \langle D \rangle^{s_{c} - \frac{1}{2}} P_{\gtrsim N} w
\|_{L_{t}^{4} L_{x}^{4} ([T,b]) }   \\
& + \| \langle D \rangle^{s_{c} - \frac{1}{2}} w \|_{L_{t}^{4} L_{x}^{4} ([T,b])}
\| w \|^{p-2}_{L_{t}^{2(p-1)} L_{x}^{2(p-1)} ([T,b])} \| P_{\gtrsim N} w \|_{L_{t}^{2(p-1)} L_{x}^{2(p-1)} ([T,b])} \\
\\
& \lesssim  \| \langle D \rangle^{1-s_{c}} I w  \|^{p-1}_{L_{t}^{2(p-1)} L_{x}^{2(p-1)} ([T,b]) }
 \frac{ \| \langle D \rangle^{1- \frac{1}{2}} I w  \|_{L_{t}^{4} L_{x}^{4}([T,b]) }}{N^{(1-s_{c})-}} \\
& + \frac{ \| \langle D \rangle^{1- \frac{1}{2}} I w  \|_{L_{t}^{4} L_{x}^{4}([T,b])} \| \langle D \rangle^{1-s_{c}} I w  \|^{p-1}_{L_{t}^{2(p-1)} L_{x}^{2(p-1)} ([T,b])  } }{N^{(1-s_{c})-}} \\
\\
&  \lesssim  \frac { \langle L_{w} \rangle^{2p(p-1)-}  A^{\frac{p(p+1)}{2}}}{N^{(1-s_{c})-}} \, .
\end{array}
\nonumber
\end{equation}
By (\ref{Eqn:DecompFLowHigh}) we write $X^{K'}_{IF(w)-F(w)} = X^{K'}_{Z_{1}} + X^{K'}_{Z_{2}} + X^{K'}_{Z_{3}}$ where

\begin{equation}
\begin{array}{l}
Z_{1} := (I-1)  F ( P_{\ll N} w ) \\
Z_{2} := (I-1) \int_{0}^{1} |
P_{\ll N} w + y P_{\gtrsim N} w|^{p-1} \,  P_{\gtrsim N} w \, dy \\
Z_{3} := (I-1)  \int_{0}^{1}  \frac{
P_{\ll N} w + y P_{\gtrsim N} w }{  _{ \quad  \overline{ P_{\ll N} w + y P_{\gtrsim N} w} \quad } } |P_{\ll N} w  + y P_{\gtrsim N} w |^{p-1} \,
\overline{P_{\gtrsim N} w} \, dy
\end{array}
\nonumber
\end{equation}

Again, we use the smoothness of $F$ ($F$ is $C^{1}$) to deal with $X^{K'}_{Z_{1}}$. We have
(using again (\ref{Eqn:BoundwlPerturb}) and (\ref{Eqn:StrNlkg}))

\begin{align*}
\begin{array}{ll}
 \| \langle  D \rangle^{s_{c} - \frac{1}{2}} X^{K'}_{Z_{1}} \|_{L_{t}^{q} L_{x}^{r} (K^{'}) }  & \lesssim
\| \langle D \rangle^{s_{c} - \frac{1}{2}} Z_{1} \|_{L_{t}^{\frac{4}{3}} L_{x}^{\frac{4}{3}} (K^{'})} \\
& \lesssim \frac{ \| \nabla F ( P_{\ll N} w ) \|_{ L_{t}^{\frac{4}{3}} L_{x}^{\frac{4}{3}} (K^{'}) } }
{N^{ \left( \frac{3}{2} - s_{c} \right)-}} \\
& \lesssim  \frac{1}{N^{ \left( \frac{3}{2} - s_{c} \right)-}} \| P_{\ll N} w \|^{p-1}_{L_{t}^{2(p-1)} L_{x}^{2(p-1)} ([T,b])} \| \nabla
P_{\ll N} w \|_{L_{t}^{4} L_{x}^{4} ([T,b])}  \\
& \lesssim \frac{1}{N^{ \left( 1 - s_{c} \right) -}}  \| \langle D \rangle^{1-s_{c}} I w \|^{p-1}_{L_{t}^{2(p-1)} L_{x}^{2(p-1)} ([T,b])}
\| \langle D \rangle^{1-\frac{1}{2}}
I w \|_{L_{t}^{4} L_{x}^{4} ([T,b]) } \\
& \lesssim \frac {  \langle L_{w} \rangle ^{2p(p-1)-}  A^{\frac{p(p+1)}{2} }}{N^{(1-s_{c})-}} \, .
\end{array}
\nonumber
\end{align*}
As for $X^{K'}_{Z_{2}}$, we have
\begin{align*}
\| \langle   &D  \rangle^{s_{c} - \frac{1}{2} } X^{K'}_{Z_{2}} \|_{L_{t}^{q} L_{x}^{r} (K^{'}) }  \lesssim \left\| \langle D \rangle^{
s_{c} - \frac{1}{2}}
Z_{2}  \right\|_{L_{t}^{\frac{4}{3}} L_{x}^{\frac{4}{3}} (K^{'})} \\
& \lesssim
 \int_{0}^{1}
\left[
\begin{array}{l}
\| \langle D \rangle^{s_{c} - \frac{1}{2}} ( |P_{\ll N} w + y P_{\gtrsim N} w|^{p-1} ) \|_{L_{t}^{\frac{4(p-1)}{3p-5}}
L_{x}^{\frac{4(p-1)}{3p-5}} ([T,b])}
\| P_{\gtrsim N} w \|_{L_{t}^{2(p-1)} L_{x}^{2(p-1)} ([T,b])} \\
+ \left\| | P_{\ll N} w + y P_{\gtrsim N} w |^{p-1} \right\|_{L_{t}^{2} L_{x}^{2} ([T,b])} \| \langle D  \rangle^{s_{c} - \frac{1}{2}}
P_{\gtrsim N} w \|_{L_{t}^{4} L_{x}^{4} ([T,b])}
\end{array}
\right] \, dy \\
& \lesssim  \frac{ \| \langle D \rangle^{\frac{1}{2}} I w \|_{L_{t}^{4} L_{x}^{4} ([T,b]) } \| \langle D \rangle^{1-s_{c}} I w
\|^{p-1}_{L_{t}^{2(p-1)} L_{x}^{2(p-1)} ([T,b])} } {N^{(1-s_{c})-}} \\
& \lesssim \frac {  \langle L_{w} \rangle ^{2p(p-1)-}  A^{\frac{p(p+1)}{2} }}{N^{(1-s_{c})-}},
\end{align*}
by using the product rule followed by a two-variable Leibnitz rule
(see Appendix B) with $f:= P_{\ll N} w$ and $g:= P_{\gtrsim N} w$, $L_{y}(f,g)=|f + y g|^{p-1}$ and $\lambda=p-2$. \\
$X^{K'}_{Z_{3}}$ is treated in a similar fashion. In fact, we get
\begin{equation*}
\| \langle  D  \rangle^{s_{c} - \frac{1}{2} } X^{K'}_{Z_{3}} \|_{L_{t}^{q} L_{x}^{r} (K^{'}) }  \lesssim \frac {  \langle L_{w} \rangle ^{2p(p-1)-}  A^{\frac{p(p+1)}{2} }}{N^{(1-s_{c})-}} \, .
\end{equation*}
Combining together, we obtain
\begin{equation}
\bar{Z}(K', I^{-1} X^{K'}_{I F(w) - F(w)} )  \lesssim \frac {  \langle L_{w} \rangle ^{2p(p-1)-}  A^{\frac{p(p+1)}{2} }}{N^{(1-s_{c})-}} \, .
\label{Eqn:barZprimew}
\end{equation}
\\
We can estimate $\bar{Z}(K', I^{-1} X_{IF(w') - F(Iw')}^{K'} )$ by  performing a similar decomposition as previously, using (\ref{Eqn:BoundwlplusonePerturb}) instead
of (\ref{Eqn:BoundwlPerturb}). We get the same bound that was found in (\ref{Eqn:barZprimew}), with $w$ substituted for $w'$.

\end{proof}

\section{Appendix B: A two-variable Leibnitz rule}
\noindent In this section we provide the proof of a two-variable Leibnitz rule.
\begin{lem}
Let $L \in C^{1} \left( \mathbb{C}^{2}, \mathbb{C} \right)$ such that $L(0,0)=0$ and such that for all $\mu \in [0,1]$ and for all
$(z_{1},z_{2},w_{1},w_{2}) \in \mathbb{C}^{4}$ we have
\begin{equation}
| L^{'} ( \mu z_{1} + (1- \mu) z_{2}, \mu w_{1} + (1-\mu) w_{2} ) | \lesssim  |z_{1}|^{\lambda} + |z_{2}|^{\lambda} +
|w_{1}|^{\lambda} + |w_{2}|^{\lambda}
\end{equation}
for some $\lambda > 0$. Then
\begin{equation}
\| L(f,g) \|_{H^{s,p}} \lesssim
\left(
\begin{array}{l}
\| f \|^{\lambda}_{L^{p_{1}}} \| f \|_{H^{s,p_{2}}} + \| g \|^{\lambda}_{L^{\tilde{p}_{1}}}  \| f
\|_{H^{s,\tilde{p}_{2}}} \\
 + \| f \|^{\lambda}_{L^{r_{1}}} \| g \|_{H^{s,r_{2}}} +  \| g \|^{\lambda}_{L^{\tilde{r}_{1}}} \| g \|_{H^{s,\tilde{r}_{2}}}
\end{array}
\right),
\label{Eqn:LeibnRule}
\end{equation}
assuming that $(p, p_{1}, p_{2}, \tilde{p}_{1}, \tilde{p}_{2}, r_{1},r_{2}, \tilde{r}_{1}, \tilde{r}_{2}) \in (1, \infty)^{9}$,
\begin{equation}
\begin{array}{l}
\frac{1}{p}  = \frac{\lambda}{p_{1}} + \frac{1}{p_{2}} = \frac{\lambda}{\tilde{p}_{1}} + \frac{1}{\tilde{p}_{2}} \; and \\
\frac{1}{p} = \frac{\lambda}{r_{1}} + \frac{1}{r_{2}} = \frac{\lambda}{\tilde{r}_{1}} + \frac{1}{\tilde{r}_{2}}.
\end{array}
\nonumber
\end{equation}
\end{lem}

\begin{proof}

The proof relies upon a simple modification of the one-variable fractional Leibnitz rules (see e.g \cite{christweins,kenigponcevega,taylor}). We recall the following inequalities (see e.g \cite{taylor} ): given $q: \mathbb{R}^{3} \rightarrow \mathbb{C}$ a function, we have
\begin{align}
\nonumber
\int_{\mathbb{R}^{3}} | P_{N_{2}} q(x) &- P_{N_{2}} q(y) | | \check{\psi_{N_{1}}} (x-y) | \, dy \lesssim \min{ \left( \frac{N_{2}}{N_{1}}, 1
\right) } M_{h} ( \tilde{P}_{N_{2}} q)(x), \\ \nonumber
\int_{\mathbb{R}^{3}} | P_{N_{2}} q(x) &-P_{N_{2}} q(y) |  | \check{ \psi_{N_{1}}} (x-y) | L(y) \, dy \\
& \lesssim \min{ \left( \frac{N_{2}}{N_{1}},1 \right)} \left(  M_{h}( \tilde{P}_{N_{2}} q )(x) M_{h} L(x)  + M_{h} ( |\tilde{P}_{N_{2}} q| L )(x)
\right)
\label{Eqn:EstHy}
\end{align}
where $(N_{1},N_{2}) \in 2^{\mathbb{N}^{*}} \in \mathbb{N}$, $H$ is a nonnegative function,  $\psi_{M}(\xi):= \psi \left( \frac{\xi}{M}\right)$ (if
$M \in 2^{\mathbb{N}}$), $\tilde{P}_{M} := P_{ \frac{M}{2} \leq \cdot \leq 2 M}$ (if $M \in 2^{\mathbb{N}^{*}}$), $\tilde{P}_{1}:= P_{\leq 2}$
and $ ( M_{h}(f) )(x):= \sup_{r >0 } \frac{1}{|B(x,r)|} \int_{B(x,r)} |f(y)| \, dy$. \\
Recall also the Paley-Littlewood inequalities (see \cite{stein})
\begin{equation*}
\| L(f,g) \|_{H^{s,p}} \lesssim \| P_{1} ( L(f,g)) \|_{L^{p}} + \left\| \left( \sum_{N_{1} \in 2^{\mathbb{N}^{*}}} N_{1}^{2s} |P_{N_{1}} (
L(f,g)) |^{2} \right)^{\frac{1}{2}} \right\|_{L^{p}}
\end{equation*}
and
\begin{equation}
\left\| \left( \sum_{N_{1} \in 2^{\mathbb{N}}} N_{1}^{2s} |\tilde{P}_{N_{1}} f|^{2} \right)^{\frac{1}{2}} \right\|_{L^{p}} \lesssim  \| f
\|_{H^{s,p}}.
\label{Eqn:Paleyf}
\end{equation}
We write
\begin{align*}
P_{N_{1}} ( L(f,g) )(x) & = \int_{\mathbb{R}^{3}} L(f(y),g(y)) \check{\psi_{N_{1}}}(x-y) \, dy   \\
& = \int_{\mathbb{R}^{3}}  \left( L(f(y),g(y)) - L(f(x),g(x)) \right)  \check{\psi_{N_{1}}}(x-y) \, dy \\
& = A_{1} + A_{2} + A_{3} + A_{4}
\end{align*}
where
\begin{align*}
A_{1} := \int_{0}^{1} \int_{\mathbb{R}^{3}} | \partial_{z} L( \mu f(y) + (1- \mu) f(x), \mu g(y) + (1- \mu) g(x) ) | |  f(y) -
f(x)) | \,  | \check{\psi_{N_{1}}}(x-y) | \, d\mu \, dy \\
A_{2} := \int_{0}^{1} \int_{\mathbb{R}^{3}} | \partial_{\bar{z}} L ( \mu f(y) + (1- \mu) f(x), \mu g(y) + (1 - \mu) g(x) ) | |
f(y) - f(x) | \,  | \check{\psi_{N_{1}}}(x-y) | \, d \mu \, dy \\
A_{3} := \int_{0}^{1} \int_{\mathbb{R}^{3}} | \partial_{w}  L ( \mu f(y) + (1- \mu) f(x), \mu g(y) + (1 - \mu) g(x) ) |  |g(y)
-g(x) |  | \check{\psi_{N_{1}}}(x-y) | \, d \mu \, dy \\
A_{4} := \int_{0}^{1} \int_{\mathbb{R}^{3}} | \partial_{\bar{w}} L ( \mu f(y) + (1- \mu) f(x), \mu g(y) + (1 - \mu) g(x) ) |
|g(y)-g(x) |  | \check{\psi_{N_{1}}}(x-y) | d \mu \, dy.
\end{align*}
Let us deal for example with $A_{1}$.
\begin{equation}
\begin{array}{ll}
\sum_{N_{1} \in 2^{\mathbb{N}^{*}}} N_{1}^{2s} A_{1}^{2}  & \lesssim
\left(
\begin{array}{l}
 \sum_{N_{1} \in 2^{\mathbb{N}^{*}}}
N_{1}^{2s} \left( \int_{\mathbb{R}^{3}} |f(y)|^{\lambda} |f(y)-f(x)|  | \check{\psi_{N_{1}}}(x-y) | \, dy \right)^{2} \\
+ \sum_{N_{1} \in 2^{\mathbb{N}^{*}}} N_{1}^{2s} |f(x)|^{2 \lambda}  \left( \int_{\mathbb{R}^{3}} |f(y)-f(x)|  | \check{\psi_{N_{1}}}(x-y) | \, dy \right)^{2} \\
+ \sum_{N_{1} \in 2^{\mathbb{N}^{*}}} N_{1}^{2s} \left( \int_{\mathbb{R}^{3}} |g(y)|^{\lambda} |f(y)-f(x)|  | \check{\psi_{N_{1}}}(x-y) | \, dy  \right)^{2} \\
 + \sum_{ N_{1} \in 2^{\mathbb{N}^{*}} } N_{1}^{2s} |g(x)|^{2 \lambda}  \left( \int_{\mathbb{R}^{3}} |f(y)-f(x)|  | \check{\psi_{N_{1}}}(x-y) |
\, dy
\right)^{2}
\end{array}
\right) \\
& \\
& \lesssim  A^{2}_{1,1} + A^{2}_{1,2} + A^{2}_{1,3} + A^{2}_{1,4}.
\end{array}
\nonumber
\end{equation}
We have
\begin{align*}
A^{2}_{1,1}  & \lesssim \sum_{N_{1} \in 2^{\mathbb{N}^{*}}} N_{1}^{2s} \left( \sum_{N_{2} \in 2^{\mathbb{N}}} \int_{\mathbb{R}^{3}}
|f(y)|^{\lambda}
|P_{N_{2}} (f)(y)- P_{N_{2}} (f)(x)| | \check{\psi_{N_{1}}}(x-y) | \, dy \right)^{2} \\
& \lesssim
\sum_{N_{1} \in 2^{\mathbb{N}^{*}}} N_{1}^{2s} \left(  \sum_{N_{2} \leq N_{1}} \int_{\mathbb{R}^{3}} |f(y)|^{\lambda} | P_{N_{2}}(f)(y) -
P_{N_{2}}(f)(x) | | \check{\psi_{N_{1}}}(x-y) | \, dy \right)^{2} \\
&+  \sum_{N_{1} \in 2^{\mathbb{N}^{*}}} N_{1}^{2s} \left(  \sum_{N_{2} \geq N_{1}} \int_{\mathbb{R}^{3}} |f(y)|^{\lambda} | P_{N_{2}}(f)(y) -
P_{N_{2}}(f)(x) | | \check{\psi_{N_{1}}}(x-y) | \, dy \right)^{2} \\
& \lesssim  A^{2}_{1,1,1} + A^{2}_{1,1,2}.
\end{align*}
But, by (\ref{Eqn:EstHy}) we have
\begin{equation}
A^{2}_{1,1,1} \lesssim
\begin{array}{l}
 ( M_{h} ( |f(x)|^{\lambda}) )^{2} \sum_{N_{1} \in 2^{\mathbb{N}^{*}} } N_{1}^{2s} \left( \sum_{N_{2} \leq N_{1}}
\frac{N_{2}}{N_{1}} M_{h} ( (\tilde{P}_{N_{2}} f)(x) )  \right)^{2} \\
+ \sum_{N_{1} \in 2^{\mathbb{N}^{*}}} N_{1}^{2s} \left( \sum_{N_{2} \leq N_{1}} \frac{N_{2}}{N_{1}} M_{h} (( |\tilde{P}_{N_{2}} f| |f|^{\lambda})
(x) ) \right)^{2}.
\end{array}
\nonumber
\end{equation}
Now, by Young's inequality  we have
\begin{align*}
\sum_{N_{1} \in 2^{\mathbb{N}^{*}}} N_{1}^{2s} \left( \sum_{N_{2} \leq N_{1}} \frac{N_{2}}{N_{1}} |a_{N_{2}}|  \right)^{2}
& = \sum_{N_{1} \in 2^{\mathbb{N}^{*}}} \left( \sum_{N_{2} \leq N_{1}} \left( \frac{N_{2}}{N_{1}} \right)^{1-s} N^{s}_{2} |a_{N_{2}}| \right)^{2} \\
& \lesssim \sum_{N_{1} \in 2^{\mathbb{N} }} N_{1}^{2s} |a_{N_{1}}|^{2}
\end{align*}
which implies that
\begin{equation}
A^{2}_{1,1,1} \lesssim   \left( M_{h} ( |f(x)|^{\lambda}) \right)^{2} \sum_{N_{1} \in 2^{\mathbb{N} }} N_{1}^{2s}  \left( M_{h} (
(\tilde{P}_{N_{1}} f)(x) ) \right)^{2} + \sum_{N_{1} \in 2^{\mathbb{N} }} N_{1}^{2s} \left( M_{h} ( (  |\tilde{P}_{N_{1}} f| |f|^{\lambda})(x) )
\right)^{2}
\label{Eqn:BoundA2111}
\end{equation}
and therefore, by Fefferman-Stein maximal inequality \cite{FefferStein}, H\"older's inequality and (\ref{Eqn:Paleyf}) we have
\begin{equation*}
\| A_{1,1,1} \|_{L^{p}} \lesssim  \| f \|^{\lambda}_{L^{p_{1}}} \| f \|_{H^{s,p_{2}}}.
\end{equation*}
Also, by (\ref{Eqn:EstHy}) we have
\begin{equation}
A^{2}_{1,1,2} \lesssim
\left(
\begin{array}{l}
\left( ( M_{h} ( |f(x)|^{\lambda}) )^{2} \sum_{N_{1} \in 2^{\mathbb{N}^{*}} } N_{1}^{2s} \left( \sum_{N_{2} > N_{1}}
 M_{h} ( (\tilde{P}_{N_{2}} f)(x) )  \right)^{2} \right) \\
+ \sum_{N_{1} \in 2^{\mathbb{N}^{*}}} N_{1}^{2s} \left( \sum_{N_{2} > N_{1}} M_{h} (( |\tilde{P}_{N_{2}} f| |f|^{\lambda}) (x) ) \right)^{2}
\end{array}
\right)
\nonumber
\end{equation}

\begin{align*}
\sum_{N_{1} \in 2^{\mathbb{N}^{*}}} N_{1}^{2s} \left( \sum_{N_{2} > N_{1}} |a_{N_{2}}|   \right)^{2} &= \sum_{N_{1} \in 2^{\mathbb{N}^{*}}}
\left( \sum_{N_{2} > N_{1}}  \left( \frac{N_{1}}{N_{2}} \right)^{s} N_{2}^{s} |a_{N_{2}}|  \right)^{2} \\
& \lesssim \sum_{N_{1} \in 2^{\mathbb{N}}} N_{1}^{2s} |a_{N_{1}}|^{2}.
\end{align*}
and therefore (\ref{Eqn:BoundA2111}) also holds if $A_{1,1,1}$ is substituted for $A_{1,1,2}$.\\
The other terms ( $A_{1,2}$, $A_{1,3}$, $A_{1,4}$ and then $A_{2}$, $A_{3}$, $A_{4}$) are treated in a similar fashion. \\
We also have $ \| P_{1}(L(f,g)) \|_{L^{p}} \lesssim  \| L(f,g) \|_{L^{p}} $. Then writing $L(f,g)= L(f,g)-L(0,0)$ and applying the fundamental
theorem of calculus, we see that (\ref{Eqn:LeibnRule}) holds if $s=0$.

\end{proof}

\subsection*{Acknowledgments}
S.K. is partially supported by NRF(Korea) grant 2010-0024017. T.R is partially supported by JSPS (Japan) grant 15K17570.

\end{document}